\documentclass[12pt,leqno]{amsart}
\usepackage[all]{xy}
\usepackage{amssymb}

\newtheorem{defn}{Definition}
\newtheorem{thm}[defn]{Theorem}

\newtheorem{cor}[defn]{Corollary}
\newtheorem{lem}[defn]{Lemma}
\newtheorem{prop}[defn]{Proposition}

\theoremstyle{remark}
\newtheorem{rem}{Remark}
\theoremstyle{remark}
\newtheorem{exam}{Example}

\numberwithin{equation}{section}
\numberwithin{defn}{section}

\renewcommand\sp{\operatorname{Spec}}
\renewcommand\sf{\operatorname{Spf}}
\newcommand\grv{{\operatorname{Gr}}(V)}
\newcommand\gr{\operatorname{Gr}}
\newcommand\lgl{\operatorname{LGl}}
\newcommand\glv{{\operatorname{Gl}}(V)}

\newcommand\gl{\operatorname{Gl}}

\renewcommand\hom{\operatorname{Hom}}

\newcommand\Detv{\operatorname{Det}_V}
\newcommand\detd{\operatorname{Det}^\ast}
\newcommand\im{\operatorname{Im}}
\newcommand\res{\operatorname{Res}}

\newcommand\limil[1]{\underset{#1}\varinjlim\,}

\newcommand\limpl[1]{\underset{#1}\varprojlim\,}
\newcommand\proj{\operatorname{Proj}}
\newcommand\aut{\operatorname{Aut}}
\newcommand\End{\operatorname{End}}
\newcommand\tr{\operatorname{Tr}}
\newcommand\Detd{\operatorname{Det}^\ast}

\newcommand\Rad{\operatorname{Rad}}

\renewcommand\o{{\mathcal O}}
\renewcommand\j{{\mathcal J}}

\renewcommand\P{{\mathcal P}}
\renewcommand\H{{\mathcal H}}
\renewcommand\L{{\mathcal L}}

\newcommand\M{{\mathcal M}}

\newcommand\Z{{\mathbb Z}}
\newcommand\A{{\mathbb A}}
\newcommand\C{{\mathbb C}}

\newcommand\w{\widehat}
\newcommand\wtilde{\widetilde}

\newcommand\beq{
                \setcounter{equation}{\value{defn}}\addtocounter{defn}1
                \begin{equation}}

\begin{document}

\title[Equations of Hurwitz schemes]{Equations of Hurwitz schemes \\ in the infinite Grassmannian}
\author[J. M. Mu\~noz Porras and F. J. Plaza Mart\'{\i}n]{J. M. Mu\~noz Porras \\  F. J. Plaza Mart\'{\i}n}
   \address{Departamento de Matem\'aticas, Universidad de Salamanca,  Plaza
   de la Merced 1-4
   \\
   37008 Salamanca. Spain.
   \\
    Tel: +34 923294460. Fax: +34 923294583}

\thanks{
 (Secondary). 
\indent This work is partially supported by the research contracts
BFM2000-1327  of DGI and  SA064/01 of JCyL. The second author is
also supported by MCYT ``Ram\'on y Cajal'' program.}
\email{jmp@@usal.es} \email{fplaza@@usal.es}

\begin{abstract}
The main result proved in the paper is the computation of the
explicit equations defining the Hurwitz schemes of coverings with
punctures as subschemes of the Sato infinite Grassmannian. As an
application, we characterize the existence of certain linear
series on a smooth curve in terms of soliton equations.
\end{abstract}


\maketitle


\section{Introduction}

It is well known that the Krichever map can be extended to the
case when the geometric data are given by a finite covering of
pointed Riemann surfaces and trivializations at the punctures.
This has been studied in works authored by M.R. Adams and M.J.
Bergvelt (\cite{AdamsBergvelt}), M.J. Bergvelt and ten A.P.E.
Kroode (\cite{BergveltnKP}) and Y. Li and M. Mulase
(\cite{LiMulasePrym,LiMulaseHitchin}). The soliton hierarchies
appearing naturally in this problem are given by the flows defined
by the Heisenberg algebras in $\w{\operatorname{sl}}(n,\C)$ (the
affine Kac-Moody algebra associated with the loop group of
$\operatorname{Sl}(n,\C)$).

The objective of this paper is to characterize the Hurwitz schemes
parametrizing finite coverings of Riemann surfaces in terms of
soliton equations satisfied by certain $\tau$-functions.

For this goal, the first main step consists of  generalizing the
``formal geometry'' developed by the authors in~\cite{MP} to the
case of the ``formal spectral cover''. Let $V$ be a finite a
separable $\C((z))$-algebra of dimension $n$ and $V_+\subset V$ a
$\C[[z]]$-subalgebra of rank $n$ over $\C[[z]]$. The properties of
the infinite Grassmannian, $\gr(V,V_+)$, and of the ``formal
spectral cover'', $\sf V_+$, are studied with the help of the
$\C((z))$-algebra structure of $V$. This construction allows us to
prove a new bilinear identity (Theorem~\ref{thm:bilinearidentity})
that depends on the algebra structure of $V$.

Consider the Hurwitz space $\H^\infty(e_1,\dots, e_r)$ (where
$e_1+\dots +e_r=n$) parametrizing geometric data $(Y,X,\pi,x,\bar
y,t_x,t_{\bar y})$, where $Y\overset{\pi}\to X$ is a finite
covering of degree $n$, $x\in X$, $\bar
y=\pi^{-1}(x)=\{e_1y_1+\dots +e_r y_r\}$ and $t_x$, $t_{\bar y}$
are formal parameters at $x$ and $\bar y$ respectively.
Theorem~\ref{thm:injectivityKricheverH} shows that the Krichever
functor defines an immersion:
$$K\colon \H^\infty(e_1,\dots, e_r)\,\longrightarrow\,
\gr(V,V_+)$$ where $V=\C((z^{1/e_1}))\times \dots\times
\C((z^{1/e_r}))$. The main result of the paper is
Theorem~\ref{thm:traceBinB}, which gives a characterization of the
image of $\H^\infty(e_1,\dots, e_r)$ in $\gr(V,V_+)$ uniquely in
terms of the piece of data $(Y,\bar y,t_{\bar y})$ of the
geometric data and the algebra structure of $V$. This
characterization allows us to prove that the Hurwtiz space is a
scheme (Theorem~\ref{thm:HisrepresentableinGRV}). Furthermore, the
$\tau$-functions of this space are explicitly characterized by an
 set of differential equations given in Theorems~\ref{thm:hurwitzequations}
and~\ref{thm:diffequationsofHinGRV}.

In the last section,  we apply the above results to study the
finite dimensional Hurwitz scheme $\H(g,0;1,\dots, 1)$, which
parametrizes finite coverings of ${\mathbb P}^1$ with a fibre of
type $(1,\dots, 1)$. One could say that the paper gives an
explicit method for constructing arbitrary finite coverings of
Riemann surfaces from a local datum (the algebra structure of $V$)
and a system of differential equations related to a soliton
hierarchy.

Another application of our results is to be found in~\cite{GMP},
where the equations defining the moduli space of curves with an
automorphism of prime order as a subscheme of the Sato
Grassmannian are given.

Along this paper we shall assume that the base field is $\C$, the
field of complex numbers.


\section{Vector-valued infinite Grassmannians}

Let $V$ be a separable and finite $\C((z))$-algebra of dimension
$n$ and let $V_+\subset V$ be a $\C[[z]]$-subalgebra of rank $n$
over $\C[[z]]$. Let us denote by $\grv$  the infinite Grassmannian
of $(V,V_+)$ constructed in~\cite{AMP} (see also~\cite{PS,SW,SS}).
It is worth recalling that \S\S2 and 3 of~\cite{AMP} are concerned
with the existence and basic properties of this Grassmannian. Let
us summarize some of them.

The infinite Grassmannian of $(V,V_+)$ is a $\C$-scheme
 whose set of rational points is:
 {\small $$\left\{\begin{gathered}
 \text{subspaces $U\subset V$ such that }U\to V/V_+
 \\
 \text{ has finite dimensional kernel and cokernel}
 \end{gathered}
 \right\}$$}
The connected components of this scheme are indexed by the
Poincar\'e-Euler characteristic of $U\to V/V_+$. The connected
component of index $m$ will be denoted by $\gr^m(V)$. Recall that
$\grv$ is equipped with the determinant bundle, whose dual has a
canonical global section $\Omega_+$.

Let us recall briefly the definition of $\Omega_+$. On the
connected component of index $0$, $\gr^0(V)$, it is the
determinant of the natural map $\L\to V/V_+$ ($\L$ being the
universal submodule). For an integer $m>0$, set $v_m\in V_+$ such
that $\dim_\C V_+/v_m V_+=m$. Then, the section $\Omega_+$ on
$\gr^m(V)$ for $m>0$ (resp. $m<0$), which will be denoted by
$\Omega^m_+$, is the determinant of the map $\L\to V/v_m V_+$
(resp. $\L\to V/v_{-m}^{-1} V_+$).

The fourth section of~\cite{AMP}  is devoted to the study of some
groups acting on the Grassmannian and uses the notions of the
formal curve, its Jacobian, etc. . Now, we shall employ the same
techniques to generalize those notions to our present setting.

\begin{exam}\label{example:couplesvv+}
The main examples of couples $(V,V_+)$ of the above type are the
following:
\begin{enumerate}
\item $V=\C((z^{1/e})),V_+=\C[[z^{1/e}]]$, where $e$ is a positive
integer.
\item $V=\C((z))\otimes_{\C}A_0,V_+=\C[[z]]\otimes_{\C}A_0$, where
$A_0$ is a finite separable $\C$-algebra.
\item $V=\C((z^{1/e_1}))\times\dots\times\C((z^{1/e_r})),
V_+=\C[[z^{1/e}]]\times\dots\times\C[[z^{1/e_r}]]$, where
$e_1,\dots,e_r$ are a positive integers.
\end{enumerate}
\end{exam}

\begin{defn}
The formal base curve associated with the couple $(V,V_+)$ is the
formal scheme $\w C:=\sf\C[[z]]$.

 The formal spectral cover associated with
the couple $(V,V_+)$ is the formal scheme:
$$\w C_V\,:=\,\sf V_+$$
\end{defn}

In the rest of this paper, it will be assumed that $\w C_V$ is a
smooth curve. Let us observe that, in general, $\w C_V$ is not
connected.

Let $V_+=V_+^1\times\dots\times V_+^r$ be the the decomposition of
$V_+$ as a product of local $\C[[z]]$-algebras. Then, the
smoothness of $\w C_V$ implies that there exist isomorphisms
$V_+^i\simeq \C[[z_i]]$ for all $i$. Further, note that the
parameters $z_i$ can be chosen such that:
$$
\begin{aligned}
\C[[z]] &\,\hookrightarrow \, V_+^i\simeq \C[[z_i]] \\
z\,&\longmapsto\, z_i^{e_i}
\end{aligned}
$$
so that one has isomorphisms $V_i\simeq\C[[z^{1/e_i}]]$. Summing
up, the  assumption of the smoothness of the  formal spectral
cover is equivalent to considering the third case of
Example~\ref{example:couplesvv+}.

{\sl Below, we shall identify $z^{1/e_i}$ with $z_i$.}

From \cite{AMP} one knows that the restricted linear group $\glv$
of the couple $(V,V_+)$, defined as a contravariant functor on the
category of $\C$-schemes, acts on $\grv$. Moreover, if $\Detv$
denotes the determinant bundle on $\grv$, then $g^*\Detv\simeq
\Detv$ for every $g\in\glv$.

The action of $\glv$ on $\grv$ induces a central extension of
functors of groups over the category of $\C$-schemes (see
\cite{AMP}, Theorem~4.3): \beq \label{eq:centralextension} 0\to
{\mathbb G}_m\to \wtilde{\glv}\to \glv\to 0 \end{equation}

Let $\underline V^*$ be the contravariant functor over the
category of $\C$-schemes with values in the category of abelian
groups defined as follows: {\small \beq \label{eq:pointsofV*}
\begin{aligned}
\underline V^*\colon
\left\{{\scriptsize\begin{gathered}\text{category of}\\
\C-\text{schemes}\end{gathered}}\right\}\,&
\rightsquigarrow\,\left\{{\scriptsize\begin{gathered}\text{category}\\\text{of groups}\end{gathered}}\right\} \\
S\,& \rightsquigarrow\,  \big(V\hat\otimes_\C H^0(S,\o_S)\big)^*
\,=\, \left\{{\scriptsize
\begin{gathered} \text{invertible elements}\\
\text{of }V\hat\otimes_\C  H^0(S,\o_S)
\end{gathered}}\right\}
\end{aligned}
\end{equation}}

Analogously as in~\cite{AMP}~\S4, one can prove that $\underline
V^*$ is representable by a formal group scheme $\Gamma_V$ and that
the connected component of the origin, $\Gamma_V^0$, decomposes
as:
$$\Gamma_V^0\,=\, \Gamma_V^-\times {\mathbb G}_m^r\times
\Gamma_V^+$$ where $\Gamma_V^-\simeq \Gamma_1^-\times \dots\times
\Gamma_r^-$, $\Gamma_V^+\simeq \Gamma_1^+\times \dots\times
\Gamma_r^+$ and $\Gamma_i=\Gamma_i^-\times{\mathbb G}_m\times
\Gamma_i^+$ is the group scheme associated with the factor
$V_i\simeq \C((z_i))$.

It follows that a point $\gamma\in\Gamma_V^0(R)$ with values in a
ring $R$ is given by a triple $(\gamma_-,\gamma_0,\gamma_+)\in
\Gamma_V^-\times {\mathbb G}_m^r\times \Gamma_V^+$: {\small
$$
\begin{aligned}
\gamma_-=(\gamma_-^{(1)},\dots,\gamma_-^{(r)}) &\text{, where }
\gamma_-^{(i)}=1+\sum_{j=-m}^{-1}a_j^{(i)}z_i^{j} \, ,\,
a_j^{(i)}\in\Rad(R) \\
\gamma_0=(\gamma_0^{(1)},\dots,\gamma_0^{(r)}) &\text{, where
$\gamma_0^{(i)}\in R$ is invertible} \\
\gamma_+=(\gamma_+^{(1)},\dots,\gamma_+^{(r)}) &\text{, where }
\gamma_+^{(i)}=1+\sum_{j\geq 1}a^{(i)}_j z^j_i \, ,\, a_j^{(i)}\in
R
\end{aligned}$$}

\begin{rem}
For the case $n=r=1$, this group was introduced in~\cite{C} and
was applied to problems related to the tame symbol.
\end{rem}

\begin{rem}
The central extension of $\glv$ of
equation~\ref{eq:centralextension} induces a central extension of
the group scheme $\Gamma_V$:
$$0\to {\mathbb G}_m\to \wtilde\Gamma_V\to \Gamma_V\to 0$$
whose restriction to $\Gamma_i=\Gamma_V^-\times{\mathbb G}_m\times
\Gamma_V^+$ is the central extension constructed in~\cite{AMP}.
Further, note that the central extension $\wtilde\Gamma_V$ gives
rise to a pairing:
$$
\begin{aligned}
\Gamma_V\times\Gamma_V \,&\longrightarrow\, {\mathbb G}_m \\
(\gamma_1,\gamma_2)& \longmapsto
\tilde\gamma_1\tilde\gamma_2\tilde\gamma_1^{-1}\tilde\gamma_2^{-1}
\end{aligned}
$$
where $\tilde\gamma_i\in\wtilde\Gamma_V$ is an element whose image
is $\gamma_i$.
\end{rem}

\begin{rem}
The algebraic version of the loop group of $\gl(n,\C)$ is the
sub-functor of groups $\underline{\lgl}(n,\C((z)))\subset \glv$
defined by: {\small $$ \underline{\lgl}(n,\C((z)))(S)\,:=\,
\left\{{\scriptsize\begin{gathered}\text{automorphisms of }V\hat \otimes H^0(S,\o_S)\\
\text{as }\C((z))\hat \otimes H^0(S,\o_S)\text{-module}
\end{gathered}}\right\}
$$}
for a $\C$-scheme $S$.

This functor is representable by $\lgl(n,\C((z)))$, which is a
formal group $\C$-scheme. It has some distinguished subgroups
$\lgl(n)^-$, $ \gl(n,\C)$ and $\lgl(n)^+$. The ``big cell'' of
$\lgl(n,\C((z)))$ is defined as: {\small $$\lgl(n)^0\,:=\,
\lgl(n)^-\cdot \gl(n,\C)\cdot \lgl(n)^+$$} and  is an open
subscheme of the connected component of the origin in
$\lgl(n,\C((z)))$ (\cite{BL} Proposition~1.11, \cite{PS}
Ch.~8.1.2, \cite{Fa}~\S2).

The theory of formal Jacobians developed in~\cite{AMP} can be
extended to the case of general formal spectral covers since the
central extension~\ref{eq:centralextension} induces a central
extension of $\lgl(n)^0$.

The natural action of $\Gamma_V$ on $V$ (by homotheties) induces a
natural action on the Grassmannian, from which one deduces a
natural immersion $\Gamma_V^0\subset \lgl(n)^0$ such that
$\Gamma_V^-\subset \lgl(n)^-$, ${\mathbb G}_m^r\subset \gl(n,k)$
and $\Gamma_V^+\subset \lgl(n)^+$. Therefore, the elements of
$\Gamma_V$ can be interpreted as matrices of size $n\times n$ with
entries in $\C((z))\subset V$
(see~\cite{AdamsBergvelt,BergveltnKP,Frenkel} for the relation of
this kind of techniques with the Heisenberg algebras).
\end{rem}

\begin{defn}
The formal Jacobian of the formal spectral cover $\w C_V$ is the
formal group scheme $\Gamma_V^-$ and will be denoted by $\j(\w
C_V)$.
\end{defn}

The above-defined Jacobian satisfies the functorial properties of
the formal Jacobian of an integral formal curve of~\cite{AMP}.
Note, moreover, that $\j(\w C_V)$ is the formal spectrum of the
ring:
$$\C\{\{x_1^{(1)},x_2^{(1)},\dots\}\}\w\otimes  \dots\w\otimes
\C\{\{x_1^{(r)},x_2^{(r)},\dots\}\}$$ where
$\C\{\{x_1^{(i)},x_2^{(i)},\dots\}\}:=\limpl{m}
\C[[x_1^{(i)},\dots, x_m^{(i)}]]$.

By the very definition of the functor $\underline{V}^*$
(see~\ref{eq:pointsofV*}) and the decomposition $\Gamma_V^-\simeq
\Gamma_1^-\times \dots\times \Gamma_r^-$, one knows that a
morphism $\w C_V\to \j(\w C_V)$  is defined by $r$ series in one
variable with coefficients in the ring:
    $$\o(\w C_V)=
    \C[[\bar z_1 ]]\times\dots\times \C[[\bar z_r]]\simeq V_+$$
where we distinguish the variables $\bar z_i\in \o(\w C_V)$ and
the variables $z_i\in \j(\w C_V)$.

Then, we define the Abel morphism of degree one:
$$
\phi_1\colon \w C_V\,\longrightarrow\, \j(\w C_V)$$ as the
morphism corresponding to the series: {\small $$
\Big((1-\frac{\bar z_1}{z_1})^{-1}, \dots,(1-\frac{\bar
z_r}{z_r})^{-1}\Big)$$}

 Equivalently, this is the
morphism induced by the ring homomorphism: {\small $$
\begin{aligned}
\C\{\{x_1^{(1)},x_2^{(1)},\dots\}\}\w\otimes  \dots\w\otimes
\C\{\{x_1^{(r)},x_2^{(r)},\dots\}\} &\to \C[[\bar
z_1]]\times\dots\times
\C[[\bar z_r]] \\
x_i^{(j)} \quad &\mapsto \quad \bar z_j^{(i)}
\end{aligned}
$$}
One checks that the triple $(\w C_V,\j(\w C_V),\phi_1)$ verifies
the Albanese property as follows from~\cite{AMP}~\S4.

\begin{rem}
The representation of the Lie algebra of $\Gamma_+^0$ as a Lie
subalgebra of $\operatorname{Lie}(\lgl(n)^0)$ is precisely the
principal Heisenberg algebra of type $\underline n=(e_1,\dots,
e_r)$. $\Gamma_V^0$, as a subgroup of $\lgl(n)^0$, is a principal
Heisenberg group of type $\underline n=(e_1,\dots, e_r)$
(\cite{AdamsBergvelt,Frenkel}).
\end{rem}

Since $\C$ has characteristic zero, we can define the exponential
map for the formal Jacobian as follows. Let $\w {\mathbb
A}_{\infty}$ be the formal group scheme $\limil{n} \w {\mathbb
A}_{n}$:
$$\w {\mathbb
A}_{\infty}\,=\,\sf \C\{\{t_1,\dots\}\}$$ endowed with the
additive group law.

Thus, the exponential map is the morphism: {\small
$$\begin{aligned} \w {\mathbb A}_{\infty}^r
\,&\overset{\exp}\longrightarrow\,\j(\w
C_V)\\
(\{a_i^{(1)}\}_{i>0},\dots,\{a_i^{(r)}\}_{i>0}) &\mapsto
\big(\exp(\sum_{i>0}\frac{a_i^{(1)}}{z_1^{i}}),\dots,\exp(\sum_{i>0}\frac{a_i^{(r)}}{
z_r^{i}})\big)
\end{aligned}
$$}
 induced by the following ring homomorphism: {\small
$$\begin{aligned} \C\{\{x_1^{(1)},\dots\}\}\w\otimes  \dots
\w\otimes  \C\{\{x_1^{(r)},\dots\}\}& \rightarrow
\C\{\{t_1^{(1)},\dots\}\}\w\otimes  \dots \w\otimes
\C\{\{t_1^{(r)},\dots\}\}\\
1\otimes\dots\otimes x_i^{(j)}\otimes \dots\otimes 1 &\mapsto
1\otimes\dots\otimes p_i(t^{(j)})\otimes \dots\otimes 1
\end{aligned}
$$}
where $t^{(j)}=(t^{(j)}_1,t^{(j)}_2,\dots)$ and $p_i(t^{(j)})$ is
the $i$-th Schur polynomial on $t^{(j)}$; that is, the coefficient
of $ z_j^{-i}$ in the series $\exp(\sum_{k>0}t^{(j)}_k z_j^{-k})$.
Obviously, the exponential map is an isomorphism of formal group
schemes.

{\it Therefore, henceforth we shall understand that the group
$\j(\w C_V)$ is the formal spectrum of the ring:
$$\C\{\{t_1^{(1)},\dots\}\}\w\otimes  \dots \w\otimes
\C\{\{t_1^{(r)},\dots\}\}$$
 and its universal element will be:
$$\prod_{i=1}^r \exp\big(\sum_{j\geq 1}\frac{t_j^{(i)}} {z_i^j}\big)$$
}

\section{$\tau$-functions and Baker-Akhiezer functions}
\label{sec:tauBAfunctions}

 Following~\S4 of~\cite{MP}, this section generalizes the notions
of $\tau$-func\-tions and Baker-Akhiezer functions) and their
properties to our situation (\cite{DJKM,K}.

Let us consider the natural action:
$$\mu\colon\Gamma_V\times\grv\,\longrightarrow\,\grv$$
given in~\cite{AMP} and let us define a Poincar\'e sheaf on
$\Gamma_V\times\grv$ as:
$$\P\,:=\, \mu^*\detd_V$$
Then, for each rational point $U\in\grv$ one has an invertible
sheaf on $\Gamma_V$:
$$\wtilde{\L}_{\tau}(U)\,:=\, \P\vert_{\Gamma_V\times\{U\}}$$
and a natural homomorphism: \beq\label{eq:morphismtau}
H^0(\Gamma_V\times\grv,\P)\,\longrightarrow\, H^0(\Gamma_V\times
\{U\},\wtilde\L_{\tau}(U))
\end{equation}

\begin{defn}
The $\tau$-section of the point $U$, $\tilde\tau_U$, is the image
of $\mu^*\Omega_+$ by the homomorphism~\ref{eq:morphismtau}. Here,
$\Omega_+$ is the canonical global section of $\Detd_V$ defined
in~\cite{AMP}.
\end{defn}

To generalize $\tau$-functions, we must restrict our definition to
the formal scheme $\j(\w C_V)=\Gamma_V^-\subset \Gamma_V$. Let us
consider the invertible sheaf on $\j(\w C_V)$:
$$\L_{\tau}(U)\,:=\, \wtilde\L_{\tau}(U)\vert_{\j(\w C_V)\times\{U\}}$$
which is trivial. A trivialization of $\L_{\tau}(U)$ can be given
by the global section:
$$\sigma_0(g)\,:=\, g\cdot \delta_U\qquad g\in\j(\w C_V)$$
where $\delta_U$ is a non-zero element in the fibre of
$\L_{\tau}(U)$ over the point $(1,U)\in \j(\w C_V)\times\{U\}$.
Then, the $\tau$-function is defined as the trivialization of the
restriction of $\tilde\tau_U$ to $\j(\w C_V)$.

\begin{defn}
The $\tau$-function of the point $U$, $\tau_U$, is the function:
$$\tau_U\in \o(\j(\w C_V))=\C\{\{t_1^{(1)},\dots\}\}\w\otimes \dots
\w\otimes  \C\{\{t_1^{(r)},\dots\}\}$$
 such that:
$$ \tau_U(g) \,=\, \frac{\Omega_+(gU)}{\sigma_0(g)}\qquad g\in \j(\w
C_V)$$
\end{defn}

\begin{rem}
Let $V_-$ be the subspace $z_1^{-1}\C[z_1^{-1}]\times \dots\times
z_r^{-1}\C[z_r^{-1}]$ and note that $V=V_-\oplus V_+$ and that
$V_-\in\grv$. Let $X\subset \grv$ denote the orbit of $V_-\subset
V$ under the action of $\Gamma_V^+$, which acts freely on $\grv$.
Therefore, the bosonization isomorphism
$B:\Omega(S)\overset{\sim}\to\o(\Gamma_V^+)$ is the isomorphism
induced by the restriction homomorphism: {\small
$$H^0(\gr^0(V),\Detd_V)\,\longrightarrow\,H^0(X,\Detd_V\vert_X)$$}
and the isomorphism
$H^0(X,\Detd_V\vert_X)\overset{\sim}\to\o(\Gamma_V^+)$ induced by
$\Omega_+$.
\end{rem}

In order to write down expressions for the Baker-Akhiezer function
analogous to the classical ones, let us observe that the
composition of the Abel morphism with the exponential map is:
{\small
$$
\xymatrix @R=5truept  @C=40truept {
 \w C_V \ar[r]^{\phi_1} & \j(\w C_V) \ar[r]^{\exp^{-1}} &
 \w{\mathbb A}_{\infty}^r
\\
 z_j\ar@{|-{>}}[r] & \phi_1(z_j) \ar@{|-{>}}[r] &  [z_j]
}
 $$}
 which maps $z_j$ to the point of $\w{\mathbb A}_{\infty}^r$ with
coordinates:
$$[z_j]\,:=\,\Big((0,\dots),\dots,(z_j,\frac{z_j^2}2,\frac{z_j^3}3,\dots),\dots,(0,\dots)\Big)$$
or, what amounts to the same, the map $\phi_1$ is induced by the
ring homomorphism: {\small $$
\begin{aligned}
\C\{\{t_1^{(1)},\dots\}\}\w\otimes \dots\w\otimes
\C\{\{t_1^{(r)},\dots\}\}&\,\to\, \C[[z_1]]\times\dots\times
\C[[z_r]]
\\
t_i^{(j)}\quad&\,\mapsto\, (0,\dots,0,\frac{z_j^i}i,0\dots,0)
\end{aligned}
$$}
It follows that there is a natural ``addition'' morphism: {\small
\beq\label{eq:additionmorphism}
 \xymatrix @R=5truept  @C=40truept
{
 \w C_V\times \j(\w C_V) \ar[r]^{\beta} & \j(\w C_V)
\\
( z_\centerdot, t) \ar@{|-{>}}[r] & t+  [z_\centerdot] }
 \end{equation}}
 where ${z_\centerdot}=(z_1,\dots,z_r)$, $t=(t^{(1)},\dots,
t^{(r)})$ and $t+  [z_\centerdot]$ denotes the point of $
\w{\mathbb
 A}_{\infty}^r$ with coordinates $(\dots,
 t^{(j)}_i+\frac{z_j^i}{i},\dots)$.

Analogously to~\cite{AMP}, one defines the Baker-Akhiezer function
of a point $U\in\grv$ (see also~\cite{Pl}). Recall that $V=\prod
V_i$.

\begin{defn}\label{defn:expressionBAj}
The $u$-th Baker-Akhiezer function of a point $U\in\grv$ is the
$V$-valued function:
    {\small $$
    \psi_{u,U}({z_\centerdot},t)\,:=\,
    \Big(\exp(-\sum_{i\geq 1}\frac{t_i^{(1)}}{z_1^i})\frac{\tau_{U_{u1}}(t+[z_1])}{\tau_U(t)}
    , \dots,
    \exp(-\sum_{i\geq 1}\frac{t_i^{(r)}}{z_r^i})\frac{\tau_{U_{ur}}(t+[z_r])}{\tau_U(t)}
    \Big)$$}
where $1\leq u\leq r$,
$U_{uv}:=(1,\dots,z_u,\dots,z_v^{-1},\dots,1)\cdot U$ and
$t+[z_v]:=(t^{(1)},\dots, t^{(v)}+[z_v],\dots, t^{(r)})$.
\end{defn}

Let $\w C_V^N$ be the formal scheme:
$$\w C_V^N\,=\,
\sf\big( (V_+^1)^{\w\otimes N}\big)\times \dots\times \sf\big(
(V_+^r)^{\w\otimes N}\big)
$$ which is the formal spectrum of the
$\C$-algebra: $$ \o(\w C_V^N)\,=\, (V_+^1)^{\w\otimes N} \w\otimes
\dots\w\otimes (V_+^r)^{\w\otimes N}
$$
 Accordingly, if we denote $(V_+^i)^{\w \otimes
N}=\C[[z_i]]^{\w\otimes N}$ by $\C[[x_1^{(i)},\dots,x_N^{(i)}]]$,
we have: {\small \beq \label{eq:decompostionOCVN} \o(\w
C_V^N)\,=\, \C[[x_1^{(1)},\dots,x_N^{(1)},\dots,
x_1^{(r)},\dots,x_N^{(r)}]]
\end{equation}}

Recall that a morphism $\w C_V^N\to \j(\w C_V)$ is determined by a
set of $r$ series $(s_1(z_1),\dots,s_r(z_r))\in \Gamma_V^-(\w
C_V^N)$. Thus, we define the Abel morphism of degree $N$ as the
morphism:
$$\phi_N\colon \w C_V^N\,\longrightarrow\, \j(\w C_V)$$
given by the series: {\small \beq \label{eq:elementAbelmorphismN}
\Big(\prod_{k=1}^N(1-\frac{x_k^{(1)}}{z_1})^{-1}, \dots,
\prod_{k=1}^N(1-\frac{x_k^{(r)}}{z_r})^{-1}\Big)
\end{equation}}

A closed point of $\grv$ is a point whose residue field is an
extension of $\C$. So, when dealing with closed points we must
consider non-finite extensions of $\C$ since local rings of $\grv$
do not need to be finitely generated. However, our arguments will
not make use of the particular structure of $\C$ and remain valid
for any such extension. Therefore, and for the sake of simplicity,
we write $\C$.

\begin{lem}\label{lem:abelNtau}
Let $U\in\gr^0(V)$ be a closed point. Let $N>0$ be an integer
number such that $V/(V_++ z_\centerdot^N U)=(0)$. Let
    $$
    \{f_i:=(f_i^{(1)}(z_1),\dots, f_i^{(r)}(z_r))  \,\vert \, 1\leq i\leq N\cdot r\}
    $$
be a basis of $V_+\cap z_\centerdot^N U$ as $\C$-vector space such
that $f^{(j)}_i\in V_+^j$.

It then holds that:
    {\tiny
    $$
    \phi_N^*\tau_U\,=\,
    \prod_{\substack{ 1\leq k<l\leq N \\
    1\leq i \leq r }}
    (x_k^{(i)}-x_l^{(i)})^{-1}
    \left| \begin{matrix}
    f_1^{(1)}(x_1^{(1)}) & \dots & f_1^{(r)}(x_1^{(r)})
    & \dots & f_1^{(1)}(x_N^{(1)}) & \dots & f_1^{(r)}(x_N^{(r)})
    \\
    \vdots & & \vdots & & \vdots & & \vdots
    \\
    f_{Nr}^{(1)}(x_1^{(1)}) & \dots & f_{Nr}^{(r)}(x_1^{(r)})
    & \dots & f_{Nr}^{(1)}(x_N^{(1)}) & \dots & f_{Nr}^{(r)}(x_N^{(r)})
    \end{matrix}\right|
    $$}
as functions of $\o(\w C_V^N)$ (up to a non-zero
constant).
\end{lem}

In the case $V=\C((z))$, this result is deeply connected with
Fay's trisecant formula ({\cite{Fay}}) and Sato's theory of
infinite Grassmann manifolds ({\cite{SS}}).

\begin{proof}
The proof consists of repeating the arguments of the proof of
Lemma~4.6 of~\cite{MP} but taking into account the
decomposition~(\ref{eq:decompostionOCVN}).

Let $A$ denote the $\C$-algebra $\o(\w C_V^N)$ and let $g$ be the
element~(\ref{eq:elementAbelmorphismN}). Then, $\phi_N^*\tau_U(x)$
is the determinant of the map: \beq \label{eq:morphismgUVV+} g
U\,\longrightarrow\, V/V_+ \end{equation}
 (from now on, we
understand that the subspaces $U,V,\dots$ have been tensorialized
by $A$).

Let $\alpha_N$ be  the following homomorphism of $A$-modules:
    {\footnotesize
    $$
    \begin{aligned}
    \alpha_N\colon V_+\,& \longrightarrow\, A^{\oplus N
    \cdot r}=  A^{\oplus r}\times\dots \times A^{\oplus r}
    \\
    (f^{(1)}(z_1),\dots, f^{(r)}(z_r)) &\mapsto \;
    \big(f^{(1)}(x^{(1)}_1),  \dots,  f^{(r)}(x^{(r)}_1)
    ,\dots, f^{(1)}(x^{(1)}_N), \dots, f^{(r)}(x^{(r)}_N)\big)
    \end{aligned}$$}
whose kernel is generated by:
    $$
    \bar g\,:=\, \big(\prod_{k=1}^N(z_1-x_k^{(1)}), \dots,
    \prod_{k=1}^N(z_r-x_k^{(r)})\big)\,\in\,  V_+
    $$

Consider the following exact sequence of complexes of $A$-modules
(written vertically):
    {\small $$ \xymatrix{ 0 \ar[r] & \bar g\cdot
    V_+ \ar[r] \ar[d]_{\beta} & V_+ \ar[r] \ar[d]_{(\pi,\alpha_N)} &
    V_+/\bar g\cdot V_+ \ar[r]
    \ar[d]_{\bar \alpha_N} & 0 \\
    0 \ar[r]&  V/z_\centerdot^N U \ar[r] & (V/z_\centerdot^N U)\oplus
    A^{\oplus N \cdot r} \ar[r] &  A^{\oplus N\cdot  r} \ar[r] & 0 }
    $$}
Observe that the complex on the l.h.s. is quasi-isomorphic to the
complex~(\ref{eq:morphismgUVV+}) and that, therefore,
$\phi_N^*\tau_U(x)$ equals the determinant of $\beta$.
Furthermore, it turns out that:
    $$
    \vert\bar\alpha_N\vert\,=\, \prod_{\substack{ 1\leq k<l\leq N \\
    1\leq i \leq r }}(x_k^{(i)}-x_l^{(i)})\,\in\, A
    $$

It remains to compute the determinant of the complex in the
middle. To this end, we consider another exact sequence of
complexes:
    {\small $$ \xymatrix{ 0 \ar[r] &  V_+\cap
    z_\centerdot^N U \ar[r] \ar[d]_{\alpha_N^U} & V_+ \ar[r]
    \ar[d]_{(\alpha_N,\pi)} & V_+/V_+\cap z_\centerdot^N U \ar[r]
    \ar[d]_{\simeq} & 0 \\
    0 \ar[r]& A^{\oplus N\cdot r} \ar[r]&  A^{\oplus N\cdot r} \oplus
    V/z_\centerdot^N U \ar[r] & V/z_\centerdot^N U \ar[r] & 0 }
    $$}
The hypothesis implies that the morphism of the complex on the
r.h.s. is an isomorphism and that its determinant belongs to
$\C^*$. In order to compute $\vert\alpha_N^U\vert$, we consider a
basis of the $A$-module $V_+\cap z_\centerdot^N U$. We may choose
a basis:
    $$
    \{f_i:=(f_i^{(1)}(z_1),\dots, f_i^{(r)}(z_r))  \,\vert \, 1\leq i\leq N\cdot r\}
    $$
as in the statement. Then, the determinant $\vert\alpha_N^U\vert$
is given by the numerator of the statement, and the conclusion
follows from the multiplicative behavior of determinants.
\end{proof}

\begin{thm}\label{thm:BAgenera}
Let $U\in\gr^0(V)$. Then:
    $$
    \psi_{u,U}({z_\centerdot},t)
    \,=\,
    (1,\dots,z_u,\dots,1) \cdot \sum_{i>0}
    (\psi_{u,U}^{(i,1)}(z_1), \dots,\psi_{u,U}^{(i,r)}(z_r))
    p_{u i,U}(t)$$
where:
    $$
    \big\{ (\psi_{u,U}^{(i,1)}(z_1), \dots,\psi_{u,U}^{(i,r)}(z_r)) \,\vert\,   i>0 \, , \, 1\leq u\leq r \big\}
    $$
is a basis of $U$ and $p_{u i,U}(t)$ are functions in $t$.
\end{thm}

\begin{proof}
Observe that the canonical morphism $\limil{N} \w C_V^N\to
\Gamma_V^-$ is a quotient by a permutation group and that $\w
C_V=\coprod \sf V_+^i$. Then, the first step consists of computing
$\psi_{u,U}\vert_{\sf V_+^i\times \w C_V^N}$ for all $N>\!>0$.
Note that the canonical morphism:
$$\xymatrix@R=50truept{\phi_{i,N}\colon \sf
V_+^i\times \w C_V^N \ar[rr]^(.6){\phi_1\times\phi_N} & &
\Gamma_V^-\times \Gamma_V^- \ar[r] &
 \Gamma_V^-}$$ is given by the series:
    {\small
    $$
    (1, \dots, (1-\frac{\bar{z_i}}{z_i})^{-1},\dots,1)\cdot g \,\in\,
    V_+\w\otimes \o(\w C_V^N)
    $$}
where $V_+^i\simeq\C[[\bar z_i]]$ and $g$ is the
element~(\ref{eq:elementAbelmorphismN}).

Accordingly, the restriction of the Baker-Akhiezer function of $U$
to the product ${\sf V_+^i\times \w C_V^N}$ is:
    $$
    g^{(i)}(z_i)^{-1}\cdot
    \frac{\phi_{i,N}^*\tau_{U_{ui}}}{\phi_{N}^*\tau_U}
    $$

In order to compute $\phi_{i,N}^*\tau_{U_{ui}}$, which coincides
with the determinant of:
    $$
    z_\centerdot^N (1,\dots,z_u,\dots,1) U \,\longrightarrow\, V/\bar g\bar g_i V_+
    $$
($\bar g_i$ being $(1, \dots, z_i - \bar{z}_i,\dots,1)$), we
proceed as in the proof of Lemma~\ref{lem:abelNtau}. We replace
$\alpha_N$ by:
    {\small $$\begin{aligned} \alpha_{i,N}(f^{(1)}(z_1), & \dots,f^{(r)}(z_r))\,:=\\ \,
    & \big(f^{(1)}(x^{(1)}_1),  \dots,  f^{(r)}(x^{(r)}_1)
    ,\dots, f^{(1)}(x^{(1)}_N), \dots, f^{(r)}(x^{(r)}_N), f^{(i)}(\bar z_i)\big)
    \end{aligned}
    $$}
which takes values in $A^{\oplus N\cdot r+1}$. We thus have that:
    $$\vert\bar\alpha_{i,N}\vert\,=\, \vert\bar\alpha_{N}\vert\cdot
    \prod_{k=1}^N(x_k^{(i)}-\bar z_i)$$

It will suffice to calculate the determinant of the following
restriction of $\alpha_{i,N}$:
    \beq\label{eq:alpharestriction}
    V_+\cap z_\centerdot^N (1,\dots,z_u,\dots,1) U \,\longrightarrow\,
    A^{\oplus N\cdot r+1}
    \end{equation}
Consider elements $g_k=(g^{(1)}_k,\dots, g^{(r)}_k)\in U$ such
that $\{g_1,\dots, g_{Nr}\}$ is a basis of $z_\centerdot^{-N}
V_+\cap U$, $\{g_1,\dots, g_{(N+1)r}\}$ is a basis of
$z_\centerdot^{-N-1}  V_+\cap U$ and $\{g_1,\dots,
g_{Nr},g_{Nr+u}\}$ is a basis of
$z_\centerdot^{-N}(1,\dots,z_u^{-1},\dots,1) V_+\cap U$. Let
$f_k:=z_\centerdot^{N} g_k$. Then,
$\{(1,\dots,z_u,\dots,1)f_k\;\vert\; k=1,\dots, Nr, Nr+u\}$ is a
basis of $V_+\cap z_\centerdot^N (1,\dots,z_u,\dots,1) U$.


The determinant of the morphism~\ref{eq:alpharestriction} is:
    {\small
    $$
    \left| \begin{array}{cccc}
     & & & \bar z_i^{\delta_{ui}} f_1^{(i)}(\bar z_i)
    \\
    M_1 & \dots & M_N & \vdots
    \\
     & & & \bar z_i^{\delta_{ui}} f_{Nr}^{(i)}(\bar z_i)
    \\
     & & & \bar z_i^{\delta_{ui}} f_{Nr+u}^{(i)}(\bar z_i)
    \end{array}\right|
    $$}
(observe that $\bar z_i^{\delta_{ui}}$ equals $\bar z_u$ if $i=u$
and $1$ otherwise) where $M_j$ is the following matrix:
    {\small $$
    M_j\,:=\,
    \left( \begin{array}{cccccccccccc}
    f_1^{(1)}(x_j^{(1)}) &  \dots &
    x_j^{(u)} f_1^{(u)}(x_j^{(u)}) & \dots &
    f_1^{(r)}(x_j^{(r)})
    \\
    \vdots & & \vdots & & \vdots
    \\
    f_{Nr}^{(1)}(x_j^{(1)}) &  \dots &
    x_j^{(u)} f_{Nr}^{(u)}(x_j^{(u)}) &  \dots &
    f_{Nr}^{(r)}(x_j^{(r)})
    \\
    f_{Nr+u}^{(1)}(x_j^{(1)}) &  \dots &
    x_j^{(u)} f_{Nr+u}^{(u)}(x_j^{(u)}) &  \dots &
    f_{Nr+u}^{(r)}(x_j^{(r)})
    \end{array}\right)
    $$}

Substituting, we have that:
    {\small $$\begin{aligned}
    \frac{\phi_{i,N}^*\tau_{U_{ui}}}{\phi_{N}^*\tau_U}
    \,=& \,
    \frac{ \bar z_i^{\delta_{ui}}}
    {\prod_{j=1}^N(\bar z_i-x^{(i)}_j)}
    \Big( f_{Nr+u}^{(i)}(\bar z_i)\prod_{j=1}^N x_j^{(u)}
    +\sum_{j=1}^{Nr} f_j^{(i)}(\bar
    z_i)p_{uj,U}(x)\Big)
    \,= \\
    \,=& \, \frac{ \bar z_i^{\delta_{ui}} \bar z_i^{-N} } {\prod_{j=1}^N(1-\frac{x^{(i)}_j}{\bar z_i})}
    \Big( f_{Nr+u}^{(i)}(\bar z_i)\prod_{j=1}^N
    x_j^{(u)}  +\sum_{j=1}^{Nr} f_j^{(i)}(\bar
    z_i)p_{uj,U}(x)\Big)\,=  \\
    \,=& \, \frac{ \bar z_i^{\delta_{ui}}  } {
    \prod_{j=1}^N(1-\frac{x^{(i)}_j}{\bar z_i})}
    \Big(g_{Nr+u}^{(i)}(\bar z_i)\prod_{j=1}^N x_j^{(u)} +\sum_{j=1}^{Nr} g_j^{(i)}(\bar
    z_i)p_{uj,U}(x)\Big)
    \end{aligned} $$}
where $p_{uj,U}(x)$ are certain polynomials which are independent
of $i$.

The statement follows now from the definitions of $\psi_{u,U}$ and
$g$ and of the variables $t$, which are a basis of the ring of
symmetric functions on the variables $x$ (the explicit expression
is $\exp(\sum_{j\geq 1} t^{(i)}_j\bar z_i^j)=\prod_{j\geq
1}(1-\frac{x_j^{(i)}}{\bar z_i})^{-1}$).
\end{proof}

Recall that, in order to define $\Omega_+^m$ and $\Omega_+^{-m}$
($m>0$), we need to choose an element $v_m$ with $\dim_\C V_+/v_m
V_+=m$. Let us set $v_m$ as follows:
\begin{itemize}
\item for $m\leq \frac12(r-n)$, let $q,p,s,t$  be  integer numbers
defined by $-m=q\cdot (n-r) +p$, $0\leq p<n-r$, $p=s\cdot r+t$,
$0\leq t< r$. Then, we set:
    {\small $$ v_m\,:=\, (z^{-1}\cdot z_\centerdot)^q  z_1^{s+1}\dots
    z_t^{s+1} z_{t+1}^s \dots z_r^s$$}
\item for $m > \frac12(r-n)$,  we set:
    {\small $$ v_m\,:=\, (z^{-1}\cdot z_\centerdot) \cdot v_{r-n-m}^{-1} $$}
\end{itemize}

\begin{thm}\label{thm:BAgeneragrmV}
Let $U\in \gr^m(V)$. It holds that:
$$
    \psi_{u,U}({z_\centerdot},t)
    \,=\,
    v_m^{-1} \cdot (1,\dots,z_u,\dots,1) \cdot \sum_{i>0}\big(\psi_{u,U}^{(i,1)}(z_1), \dots,
    \psi_{u,U}^{(i,r)}(z_r)\big)p_{u i,U}(t)$$
where:
    $$
    \big\{ (\psi_{u,U}^{(i,1)}(z_1), \dots,\psi_{u,U}^{(i,r)}(z_r)) \,\vert\,   i>0 \, , \, 1\leq u\leq r \big\}
    $$
is a basis of $U$ and $p_{u i,U}(t)$ are functions in $t$.
\end{thm}

\begin{proof}
The commutative diagram:
    {\small $$\xymatrix{
    v_m^{-1}U \ar[r] \ar[d]_{\simeq} & V/V_+ \ar[d]_{\simeq} \\
    U \ar[r] & V/v_m V_+ }$$
shows that
$\tau_{v_m^{-1}U}(g)=\tau_U(g)$. We thus have:
    $$\begin{aligned}
    \psi_{u,U}(g,z)\,=& \, \psi_{u,v_m^{-1}U}(g,z)\,= \\
    =&\,  (1,\dots,z_u,\dots,1)\cdot
    \sum_{i>0}\big(\psi_{u,v_m^{-1}U}^{(i,1)}(z_1), \dots,
    \psi_{u,v_m^{-1}U}^{(i,r)}(z_r)\big)p_{u i,v_m^{-1}U}(t)\,= \\
    =&\, v_m^{-1} (1,\dots,z_u,\dots,1)\cdot
    \sum_{i>0}\big(\psi_{u,U}^{(i,1)}(z_1), \dots,
    \psi_{u,U}^{(i,r)}(z_r)\big)p_{u i,U}(t)
    \end{aligned}$$}
where $\psi^{(i,j)}_{u,U}(z_j):=v_m\cdot
\psi^{(i,j)}_{u,v_m^{-1}U}(z_j)$. In particular, observe that:
    {\small $$
    \big\{ (\psi_{u,U}^{(i,1)}(z_1), \dots,\psi_{u,U}^{(i,r)}(z_r)) \,\vert\,   i>0 \, , \, 1\leq u\leq r \big\}
    $$}
 is a basis of $U$.
\end{proof}


The above results allow us to prove a generalization of the
bilinear identities of the KP-hierarchy.

Since $V$ is a finite separable $\C((z))$-algebra, it carries the
metric of the trace $\tr: V\times V \to \C((z))$, which is
non-degenerate. Therefore, $V$ can be endowed with the
non-degenerate pairing:
$$
\begin{aligned}
T_2\colon V\times V &\,\longrightarrow\, \C \\
(a,b)\,&\,\longmapsto \, \res_{z=0}(\tr(a,b))dz
\end{aligned}
$$

\begin{lem}
The pairing $T_2$ gives rise to an isomorphism of $\C$-schemes:
$$
\begin{aligned}
R\colon \grv &\,\longrightarrow\, \grv \\
U&\,\longmapsto \, U^{\perp}
\end{aligned}
$$
where $U^{\perp}$ denotes the orthogonal of $U$ w.r.t. $T_2$.
\end{lem}

\begin{proof}
The proof reduces easily  to the $r=1$ case; that is,
$V=\C((z^{1/e}))$. Then, the metric given by the trace is:
    $$
    \tr(z^{i/e}, z^{j/e})\,=\,
    \begin{cases} e  & \text{ if }i+j=0
        \\ e z & \text{ if }i+j=e
        \\  0 & \text{ otherwise}
    \end{cases} \qquad \text{ for }0\leq i,j<e
    $$
Taking into account that $V$ is a $\C((z))$-algebra, one has that:
$$T_2(z^{i/e}, z^{j/e})\,=\,\begin{cases} e  & \text{ if
}i+j=-e \\ 0 & \text{ otherwise} \end{cases} \qquad \text{ for }
i,j\in\Z$$

A straightforward calculation shows that $U^\perp$ belongs to
$\grv$ for any $U\in\grv$; that is, we have obtained the morphism
$R$ (\cite{MP},~\S5).
\end{proof}

It is worth pointing out the following identities: {\small
$$\begin{aligned}
 R(\gr^m(V))\,&=\,\gr^{r-n-m}(V) \\
R^*\Detv\,&\simeq\,\Detv
\end{aligned}$$}
and $ (g\cdot U)^\perp =g^{-1}\cdot U^\perp$  for $U\in\grv$ and
$g\in\j(\w C_V)$. Further, for $m\neq\frac12(r-n)$ and
$U\in\gr^m(V)$, it holds that:
 {\small $$\begin{aligned}
R^*\Omega_+^m\,&=\,\Omega_+^{r-n-m}
\\
\tau_{U^{\perp}}(g)\,& =\, \tau_U(g^{-1})
\end{aligned}$$}

\begin{rem}
The latter two identities also hold for $m=\frac12(r-n)$ whenever
one can take $v_m$ such that $v_m^2=z^{-1}\cdot z_\centerdot$
(e.g. when $e_i$ is odd for all $i$). Since this is not posible in
general, we will omit this case. However, although with different
explicit expressions, our techniques can be applied to it.
\end{rem}

\begin{defn}
The $u$-th adjoint Baker-Akhiezer function of a point $U\in\grv$
is defined by:
$$\psi_{u,U}^*({z_\centerdot},t)\,:=\, \psi_{u,U^\perp}({z_\centerdot},-t)$$
\end{defn}

Note that one has the following identity:
    {\small $$
    \psi_{u,U}^{*(j)}(z_j,t)\,=\,
    \exp(\sum_{i\geq 1}\frac{t_i^{(j)}}{z_j^i}) \frac{\tau_{U_ju}(t-[z_j])}{\tau_U(t)} $$}

\begin{thm}[Bilinear Identity]\label{thm:bilinearidentity}
Let $U,U'\in \gr^m(V)$ ($m\neq\frac12(r-n)$) be two rational
points of the same index.

Then, $U=U'$ if and only if the following conditions hold:
    {\small
    $$
    T_2\big(
    \frac{z_\centerdot\cdot \psi_{u,U}({z_\centerdot},t)}{(1,\dots,z_u,\dots,1)},
    \frac{\psi_{v,U'}^*({z_\centerdot},t')}{z (1,\dots,z_v,\dots,1)}\big)\,=\,0
    \qquad 1\leq u,v\leq r
    $$}
\end{thm}

\begin{proof}
Since $T_2$ is non-degenerate we know that a vector $w \in V$ lies
in $U'\in\gr^m(V)$ if and only if:
    $$T_2\big(w,
    \frac{v_{r-n-m} \psi_{v,U'}^*({z_\centerdot},t')}
    {(1,\dots,z_v,\dots,1)}\big)\,=\,0
    \qquad 1\leq v\leq r
    $$
Recalling that $v_m v_{r-n-m}=z^{-1} z_\centerdot$ and the
properties of the trace, the conclusion follows from
Theorem~\ref{thm:BAgeneragrmV}.
\end{proof}

\begin{cor}\label{cor:nKP}
Let $U,U'$ be two rational points of $U,U'\in \gr^m(V)$
($m\neq\frac12(r-n)$). Then, $U=U'$ if and only if:
    $$
    \sum_{i=1}^r
    \res_{z=0}\Big(
    \sum_{j=1}^{e_i}\big(\xi_i^{j}z^{1/e_i}\big)^{1-\delta_{iu}-\delta_{iv}}
    \psi_{u,U}^{(i)}(\xi_i^j z^{1/e_i} ,t)\psi_{v,U'}^{*(i)}(\xi_i^j
    z^{1/e_i},t')\Big)\frac{dz}z\,=\,0
    $$
for all $1\leq u,v\leq r$ ($\xi_i$ is a primitive $e_i$-th root of
$1$ in $\C$).
\end{cor}

\begin{proof}
It suffices to make explicit the condition of the previous
theorem. The very definition of the metric $T_2$ yields:
    {\small
    $$
    \sum_{i=1}^r
    \res_{z=0}\Big( \tr\big( z_i^{1-\delta_{iu}-\delta_{iv}} \psi_{u,U}^{(i)}(z_i,t)
    \psi_{v,U'}^{*(i)}(z_i,t')\big)\Big)\frac{dz}z$$}
and the claim follows
since the trace map of $V$ as a $\C((z))$-algebra is given by:
{\small $$
\begin{aligned}
\tr\colon V=V_1\times\dots\times V_r &\longrightarrow \, \C((z))
\\
(f_1(z_1),\dots,f_r(z_r))&\longmapsto \sum_{j=1}^{e_1}
f_1(\xi_1^jz^{1/e_1})+\dots+ \sum_{j=1}^{e_r}
f_r(\xi_r^jz^{1/e_r})
\end{aligned}$$}
\end{proof}

 This set of equations is equivalent to a set of differential equations
for the $\tau$-functions.

\begin{defn}
Let $\underline n$ denote the partition of $n$ given by
$\{e_1,\dots, e_r\}$. The $\underline n$-KP hierarchy is the
following set of equations:
    $$
    \sum_{i=1}^r
    \res_{z=0}\Big(
    \sum_{j=1}^{e_i}\big(\xi_i^{j}z^{1/e_i}\big)^{1-\delta_{iu}-\delta_{iv}}
    \psi_{u,U}^{(i)}(\xi_i^j z^{1/e_i} ,t)\psi_{v,U}^{*(i)}(\xi_i^j
    z^{1/e_i},t')\Big)\frac{dz}z\,=\,0
    $$
for all $1\leq u,v\leq r$ ($\xi_i$ is a primitive $e_i$-th root of
$1$ in $\C$).
\end{defn}

It would be interesting to compare the other hierarchies with
those given in~\cite{BergveltnKP}, which are expressed in terms of
pseudodifferential operators and representation theory of infinite
dimensional Lie algebras.


\subsection{Subschemes of the Grassmannian}

For each subset $i_\centerdot=\{i_1,\dots,i_s\}$ of
$\{1,\dots,r\}$, let us denote by $V_{i_\centerdot}$ the vector
space $\prod_{j\in i_\centerdot}  V_j$,
$(V_{i_\centerdot})_+=\prod_{j\in i_\centerdot}(V_j)_+$, by
$V^{i_\centerdot}$ the vector space $\prod_{j\notin i_\centerdot}
V_j$ and by $(V^{i_\centerdot})_+=\prod_{j\notin i_\centerdot}
V^j_+$. One can consider the following morphism:
\beq\label{eq:producttogrv}
\begin{aligned}
\gr(V_{i_\centerdot})\times\gr(V^{i_\centerdot})
&\,\overset{\gamma_{i_\centerdot}}\longrightarrow\, \grv
\\
(W,W') &\longmapsto W\times W'\subset V
\end{aligned}
\end{equation}

\begin{defn}\label{defn:decomposablegrv}
A subspace $U$ is said ``decomposable'' if it lies in the image of
$\gamma_{i_\centerdot}$ for some ${i_\centerdot}$. That is, there
exists a subset ${i_\centerdot}$ and subspaces
$W\in\gr(V_{i_\centerdot}),W'\in\gr(V^{i_\centerdot})$ such that
$U=W\times W'$.

The ``decomposable Grassmannian'' of $V$ is the subscheme of
$\grv$ whose points are the decomposable subspaces, that is:
$$\gr^{\text{dec}}(V)\,=\, \underset{{i_\centerdot}\subsetneq \{1,\dots,r\}}\bigcup \im\gamma_{i_\centerdot}$$
($\im\gamma_{i_\centerdot}$ being the scheme-theoretic image of
$\gamma_{i_\centerdot}$).
\end{defn}

\begin{prop}\label{prop:grdec}
The morphism~\ref{eq:producttogrv} is a closed immersion for any
$i_\centerdot$. In particular, $\gr^{\text{dec}}(V)$ is a closed
subscheme of $\grv$.
\end{prop}

\begin{proof}
The map is clearly injective. If $U$ denotes an $S$-valued point
of $\grv(S)$, one has to show that the subset of $S$ of those
$s\in S$ such that $U_s$ decomposes as a product of subspaces is a
closed subscheme of $S$.

Let $p_{i_\centerdot}$ (resp. $p^{i_\centerdot}$) denote the
projection $V\to V_{i_\centerdot}$ (resp. $V\to
V^{i_\centerdot}$). Then, there is a natural injective morphism:
$$U\,\longrightarrow\, p_{i_\centerdot}(U)\times p^{i_\centerdot}(U)\subseteq V_{i_\centerdot}\times
V^{i_\centerdot}=V$$ The desired subset consists exactly of those
points where
 $p_{i_\centerdot}(U)\times p^{i_\centerdot}(U)\subseteq U$, which is a
closed subscheme.
\end{proof}

One can compute explicitly the equations of all these closed
subschemes of $\grv$:

\begin{thm}\label{thm:characterizationgrred}
Let $U$ be a closed point of $\gr^m(V)$ ($m\neq\frac12(r-n)$). It
holds that:
\begin{enumerate}
\item
$U\in \im\gamma_{i_\centerdot}$ if and only if its Baker-Akhiezer
function satisfies the following equations:
    $$\res_{z=0}
    \Big(\sum_{i\in i_\centerdot}
    \tr\big(z_i^{1-\delta_{iu}-\delta_{iv}}
    \psi_{u,U}^{(i)}(z_i,t)\psi_{v,U}^{*(i)}(z_i,t')\big)\Big)\frac{dz}z
    \,=\,0$$
for all $1\leq u,v\leq r$.
\item
$U\in \gr^{\text{dec}}(V)$ if and only if its Baker-Akhiezer
function satisfies the following equations:
    $$
    \prod_{{i_\centerdot}\subsetneq \{1,\dots,r\}} \res_{z=0}
    \Big(\sum_{i\in
    i_\centerdot}\tr\big(z_i^{1-\delta_{iu}-\delta_{iv}}
    \psi_{u,U}^{(i)}(z_i,t)\psi_{v,U}^{*(i)}(z_i,t')
    \big)\Big)\frac{dz}z\,=\,0
    $$
for all $1\leq u,v\leq r$.
\end{enumerate}
\end{thm}

\begin{proof}
This follows from Theorem~\ref{thm:bilinearidentity} and the
Bilinear Identity for the KP hierachy.
\end{proof}


\section{Algebro-geometric points of $\grv$}

The goal of this section consists of giving an explicit
characterization of some points of $\grv$ defined by geometric
data (see \cite{AdamsBergvelt,LiMulasePrym,LiMulaseHitchin,MP}).
More precisely, we wish to define a subscheme of $\grv$
representing the Hurwitz functor of pointed coverings with formal
parameters at the marked points.

In fact, we restrict ourselves to the following type of coverings.
Let $\pi:Y\to X$ be a finite morphism between proper curves over
$\C$ where $Y$ is reduced and $X$ integral. Let $x\in X$ be a
smooth point. Define $A:=H^0(X-x,\o_X)$,
$B:=H^0(Y-\pi^{-1}(x),\o_Y)$, $\Sigma_X$ the function field of $A$
and $\Sigma_Y$ the total quotient ring of $B$, and let $\tr$
denote the trace of $\Sigma_Y$ as a finite $\Sigma_X$-algebra.

{\it The triple $(Y,X,x)$ is said to have the property $(\ast)$
when $\tr(B)\subseteq A$}.

Let us observe that every covering $\pi:Y\to X$ with either $X$
smooth or $\pi$ flat has the property $(\ast)$.

 From now on, we set the numerical invariants $n$ and
$\underline n=\{e_1,\dots , e_r\}$ (with $e_1+\dots+e_r=n$), which
define the $\C((z))$-algebra $V$.

\begin{defn}
The Hurwitz functor $\underline \H^{\infty}$ of pointed coverings
of curves of degree $n$ with a fibre of type $(e_1,\dots, e_r)$
and a formal parameter along the fibre is the contravariant
functor on the category of $\C$-schemes:
$$
\underline \H^{\infty}\colon
\left\{{\scriptsize\begin{gathered}\text{category of}\\
\C-\text{schemes}\end{gathered}}\right\}\,
\rightsquigarrow\,\left\{{\scriptsize\begin{gathered}\text{category}\\\text{of
sets}\end{gathered}}\right\}
$$
that associates with a $\C$-scheme $S$ the set of equivalence
classes of data $\{Y,X,\pi,x,\bar y,t_x,t_{\bar y}\}$ where:
\begin{enumerate}
\item
$p_Y:Y\to S$ is a proper and flat morphism whose fibres are
geometrically reduced curves.
\item
$p_X:X\to S$ is a proper and flat morphism whose fibres are
geometrically integral curves.
\item
$\pi:Y\to X$ is a finite morphism of $S$-schemes of degree $n$
such that its fibres over closed points $s\in S$ have the property
$(\ast)$.
\item
$x:S\to X$ is a rational $S$-point such that the divisor $x(s)$ is
a smooth point of $X_s:=p_X^{-1}(s)$ for all closed points $s\in
S$.
\item
$\bar y=\{y_1,\dots, y_r\}$ is a set of disjoint smooth sections
of $p_Y$ such that the Cartier divisor $\pi^{-1}(x(S))$ is
$e_1y_1(S)+\dots+ e_r y_r(S)$.
\item
For all closed point $s\in S$ and each irreducible component of
the fibre $Y_s$, there is at least one point $y_j(s)$ lying on
that component.
\item
$t_x$ is a formal parameter along $x(S)$:
$$t_x\colon \w \o_{X,x(S)}\,\overset{\sim}\to\, \o_S[[z]]$$
\item
$t_{\bar y}=\{t_{y_1},\dots, t_{y_r}\}$ are formal parameters
along $y_1(S),\dots, y_r(S)$ such that
$\pi^*(t_x)_{y_j(S)}=t_{y_j}^{e_j}$.
\item
$\{Y,X,\pi,x,\bar y,t_x,t_{\bar y}\}$ and $\{Y',X',\pi',x',\bar
y',t_{x'},t_{\bar y'}\}$ are said to be equivalent when there is a
commutative diagram of  $S$-schemes: {\small
$$\xymatrix{ Y \ar[d]_{\pi} \ar[r]^{\sim} & Y' \ar[d]^{\pi'} \\
X \ar[r]^{\sim} & X'}$$} compatible with all the data.
\end{enumerate}
\end{defn}

Now, a Krichever morphism can be defined for this functor as the
morphism of functors: \beq\label{eq:kricheverfunctorK} K\colon
\underline \H^{\infty}(e_1,\dots,e_r)\,\longrightarrow
\,\underline \gr(V) \end{equation}
 given by:
$$K(Y,X,\pi,x,\bar y,t_x,t_{\bar y})\,=\,
t_{\bar y}\Big(\limil{i} (p_Y)_*\o_Y(i
\pi^{-1}(x))\Big)\,\subset\,
 V\hat\otimes_\C \o_S$$ where $\o_Y(i
\pi^{-1}(x))$ is the sheaf associated with the Cartier divisor $i
\pi^{-1}(x)$ and  $t_{\bar y}$ is understood as the isomorphism
induced by:
$$\w\o_{Y,y_1(S)}\times\dots \times \w\o_{Y,y_r(S)}\,\simeq\,
\o_S[[z_1]]\times\dots\times\o_S[[z_r]] \,\simeq\,
V_+\w\otimes_\C\o_S$$

Note that for a closed point $(Y,X,\pi,x,\bar y,t_x,t_{\bar y})\in
\underline \H^{\infty}(e_1,\dots,e_r)$ these definitions yield:
$$K(Y,X,\pi,x,\bar y,t_x,t_{\bar y})\,=\,
t_{\bar y}\big(H^0(Y-\pi^{-1}(x), \o_Y)\big)\,\subset\, V$$

Let $\M^{\infty}(r)$ be the moduli scheme representing the classes
of  sets of data $(Y; y_1,\dots, y_r; t_1,\dots, t_r)$ of
geometrically reduced curves with $r$ marked pairwise distinct
smooth points $\{y_1,\dots,y_r\}$ and formal parameters
$\{t_1,\dots , t_r\}$ at these points and such that any
irreducible component contains at least one of the marked points.
Following the arguments of \cite{MP} for $\M^{\infty}(1)$, we can
prove that the Krichever morphism induces a closed immersion:
$$
\begin{aligned}
\M^{\infty}(r) \,&\overset{K}\hookrightarrow\, \grv
\\
(Y; y_1,\dots, y_r; t_1,\dots, t_r)& \mapsto t_{\bar y}\big(
H^0(Y-\{y_1,\dots, y_r\},\o_Y)\big)\subset V
\end{aligned}
$$
whose image is characterized by the following:
\begin{thm}\label{thm:minftyr}
A point $U\in\grv(S)$ lies in $K(\M^{\infty}(r))$ if and only if
$U\cdot U\subseteq U$ and $\o_S\subseteq U$, where $\cdot$ denotes
the product of $V$.
\end{thm}

\begin{proof}
The direct proof is trivial. Let us prove the converse. Consider
the filtration of $V=V_1\times\dots \times
V_r=\C((z_1))\times\dots\times \C((z_r))$ defined by:
$$ \dots \, \subset \, V(m-1) \,\subset \, V(m) \,\subset \, V(m+1)
\,\subset\, \dots
$$
where $V(m):= z^{-m}V_+$.

Then, any point $U\in\grv(R)$ ($R$ being a $\C$-algebra) carries a
natural filtration $\{U(m):=U\cap V(m)\}_{m\geq 0}$. Let us denote
by ${\mathcal U}$ the corresponding graded $R$-module. If $U$
satisfies $U\cdot U\subseteq U$ and $R\subseteq U$, then
${\mathcal U}$ is also a graded $R$-algebra.

It is easy to check that $Y=\proj{\mathcal U}$ is an algebraic
curve over $R$. Observe that the filtrations induced by
$z_1,\dots,z_r$ give rise to pairwise disjoint sections of $Y$
(smooth and of degree $1$). The other geometric data are
constructed using the same arguments as in the proof of
Theorem~6.4 of~\cite{MP}.
\end{proof}

Let us denote the trace map of the separable $\C((z))$-algebra $V$
by:
$$\tr\colon V\,\longrightarrow\, \C((z))$$
which is a $\C((z))$-linear map. For a point $U\in\grv$, let us
denote by $\tr(U)\subseteq\C((z))$ the image of $U$ under the
trace map.

Note that for a point ${\mathcal Y}:= (Y,X,\pi,x,\bar
y,t_x,t_{\bar y})$ of $\underline\H^\infty(e_1,\dots,e_r)$, we
have a commutative diagram: {\small $$\xymatrix{
B \ar@{^(->}[r] & V \\
A \ar@{^(->}[u] \ar@{^(->}[r] & \C((z)) \ar@{^(->}[u] }$$} where
 $B:=K({\mathcal Y})= t_{\bar
y}(H^0(Y-\pi^{-1}(x),\o_Y))$  and $A:=t_x(H^0(X-x,\o_X))$. Then,
the restriction of the trace of $V$ as a $\C((z))$-algebra to $B$
coincides with the restriction of the trace of $\Sigma_Y$ as a
$\Sigma_X$-algebra to $B$. Therefore, both trace maps will be
denoted by $\tr$.

\begin{lem}\label{lem:traceinBimpliesinGrV}
Let ${\mathcal B}$ be a $S$-valued point of $\grv$ such that
$\tr({\mathcal B})\subseteq {\mathcal B}$.

It holds that:
    $$
    \tr({\mathcal B})\,\in\,\gr(\C((z)))(S)
    $$
and $\tr({\mathcal B})={\mathcal B}\cap \o_S((z))$.
\end{lem}

\begin{proof}
By the local nature of the hypotheses, we may assume that $S$ is
affine, $S=\sp R$. Since ${\mathcal B}\in\grv$, there exists $m$
such that $\w V_S/({\mathcal B}+z^m\cdot \w V^+_S)=(0)$ and that
${\mathcal B}\cap z^m\cdot \w V^+_S$ is a free $R$-module of
finite rank.

Observe that $\tr({\mathcal B})$ is quasicoherent and that
$\tr({\mathcal B})_s\subset \w V_s$ for all closed point $s\in S$.
So, in order to show that $\tr({\mathcal B})\in\gr(\C((z)))$ it
suffices to check that $R((z))/(\tr({\mathcal B})+z^m\cdot
R[[z]])=(0)$ and that $\tr({\mathcal B})\cap z^m\cdot R[[z]]$ is
free of finite rank (see~\cite{AMP}).

For the first claim, note that the trace gives rise to a
surjection:
    $$
    \w V_S/({\mathcal B}+z^m\cdot \w V^+_S)\,
    \overset{\frac1n\tr}\longrightarrow\,
    R((z))/(\tr({\mathcal B})+z^m\cdot R[[z]])
    $$
The second claim follows from the fact that the composition:
    $$
    \tr({\mathcal B})\cap z^m\cdot R[[z]]
    \,\hookrightarrow\,
    {\mathcal B}\cap z^m\cdot \w V^+_S
    \,\overset{\frac1n\tr}\longrightarrow\,
    \tr({\mathcal B})\cap z^m\cdot R[[z]]
    $$
is the identity map because $\tr(z^m\cdot \w V^+_S)=z^m \tr( \w
V^+_S)=z^m\cdot R[[z]]$.

The second part of the statement follows easily from the
$R((z))$-linearity of $\tr$ and from $\tr({\mathcal B})\subseteq
{\mathcal B}$.
\end{proof}

\begin{lem}\label{lem:traceBequalArelative}
Let ${\mathcal Y}:= (Y,X,\pi,x,\bar y,t_x,t_{\bar y})$ be an
$S$-valued point of $\underline\H^\infty(e_1,\dots,e_r)$.

It holds that:
    {\small $$ K(X,x,t_x)\,=\, K({\mathcal Y}) \cap
    \o_S((z))\,=\,  \tr(K({\mathcal Y}))$$}
\end{lem}

\begin{proof}
As in the previous lemma we may assume that $S$ is affine, $S=\sp
R$. For the point ${\mathcal Y}$, define ${\mathcal
B}:=K({\mathcal Y})= t_{\bar y}(H^0(Y-\pi^{-1}(x),\o_Y))$  and
${\mathcal A}:=K(X,x,t_x)= t_x(H^0(X-x,\o_X))$.

From the commutative diagram: {\small
$$\xymatrix{
{\mathcal B} \ar@{^(->}[r] & \w V_S \\
{\mathcal A} \ar@{^(->}[u] \ar@{^(->}[r] & R((z)) \ar@{^(->}[u]
}$$} one has that ${\mathcal A}\subseteq {\mathcal B}\cap R((z))$.
The inclusion ${\mathcal B}\cap R((z))\subseteq{\mathcal  B}$
implies that $\tr({\mathcal B}\cap R((z)))\subseteq \tr({\mathcal
B})$. Having in mind that $\tr$ is $R((z))$-linear, one concludes
that ${\mathcal B}\cap R((z))\subseteq \tr({\mathcal B})$. Summing
up, we have proved that:
$${\mathcal A}\,\subseteq\, {\mathcal B}\cap R((z))\,\subseteq \,
\tr({\mathcal B})$$

Since ${\mathcal B}$ is a finite ${\mathcal A}$-module, so
$\tr({\mathcal B})$ does. Bearing in mind the compatibility of the
trace w.r.t. base changes, the above inclusion implies that
${\mathcal A}_s\subseteq \tr({\mathcal B})_s$ for all closed
points $s\in S$. Now, recalling that the data ${\mathcal Y}$
satisfies the property $(\ast)$ at closed points, one obtains that
${\mathcal A}_s\subseteq \tr({\mathcal B})_s$ for all $s$ and
that, therefore:
    $$
    {\mathcal A}\,=
    \, {\mathcal B}\cap R((z))\,
    = \, \tr({\mathcal B})
    $$
\end{proof}

\begin{thm}\label{thm:injectivityKricheverH}
The Krichever morphism~\ref{eq:kricheverfunctorK} is injective.
\end{thm}

\begin{proof}
We will keep the same notations as in the previous lemma. It
suffices to show that the geometric data ${\mathcal Y}:=
(Y,X,\pi,x,\bar y,t_x,t_{\bar y})$ can be recovered from the point
${\mathcal B}=K({\mathcal Y})\in\grv$. Observe that $\mathcal B$
determines uniquely the data $(Y,\bar y,t_{\bar y})$.

Lemma~\ref{lem:traceBequalArelative} shows that:
$${\mathcal A}\,=\, \tr({\mathcal B})\,=\,{\mathcal B}\cap\o_S((z))\, \in\,
\gr(\C((z)))(S)$$ and note that this is an $\o_S$-subalgebra of
$\o_S((z))$ because ${\mathcal B}$ is an $\o_S$-subalgebra of
$V\w\otimes_\C\o_S$.

Theorem~\ref{thm:minftyr} shows that ${\mathcal A}=\tr({\mathcal
B})$ corresponds to the point $(X,x,t_{x})$ of $\M^\infty(1)(S)$.
Since the inclusion $\tr({\mathcal B})\subseteq {\mathcal B}$ is
compatible with the filtrations induced by those of $\C((z))$ and
$V$ respectively, it gives rise to a morphism $\pi:Y\to X$, which
is finite and such that $\pi^{-1}(x)=\bar y$ and
$(\pi^*t_{x})_{y_j}=t_{y_j}^{e_j}$.
\end{proof}

All the above results allow us to give characterizations of the
functor $\underline \H^{\infty}(e_1,\dots,e_r)$ as a subset of
$\M^{\infty}(r)\subset\underline\gr(V)$.

\begin{cor}\label{cor:UcapC((z))}
Let ${\mathcal B}$ be an $S$-valued point of $\underline
\M^{\infty}(r)$. The point ${\mathcal B}$ belongs to $\underline
\H^{\infty}(e_1,\dots,e_r)(S)$ if and only if:
$${\mathcal B}\cap \o_S((z)) \,\in\,\underline\gr(\C((z)))(S)$$
\end{cor}

\begin{proof}
The first part of the proof follows from
Lemma~\ref{lem:traceBequalArelative}. The converse is a
consequence of the proof of
Theorem~\ref{thm:injectivityKricheverH}.
\end{proof}

\begin{thm}\label{thm:traceBinB}
Let ${\mathcal B}\in\M^{\infty}(r)\subset \grv$ be an $S$-valued
point. The following conditions are equivalent:
\begin{enumerate}
\item
${\mathcal B}\in \underline\H^{\infty}(e_1,\dots,e_r)(S)$.
\item
$\tr({\mathcal B})\subseteq {\mathcal B}$.
\end{enumerate}
\end{thm}

\begin{proof}
(1)$\implies$(2) is a consequence of
Lemma~\ref{lem:traceBequalArelative}.

(2)$\implies$(1) follows from Lemma~\ref{lem:traceinBimpliesinGrV}
and Corollary~\ref{cor:UcapC((z))}.
\end{proof}

\begin{rem}
Our approach to the Hurwitz functor is closely related to the Li
and Mulase study of  the category of morphisms of algebraic curves
(\cite{LiMulasePrym}). Those authors study the equivalence of some
geometrical data (essentially a pointed covering with parameters
and a vector bundle upstairs) and certain triples $(A,B,W)$ of
points of the Grassmannian (where $A$ is related to the curve
downstairs, $B$ to the curve upstairs and $W$ to the vector
bundle). From their point of view, we are restricting ourselves to
the case when the vector bundle is the sheaf of algebraic
functions on the curve upstairs and the covering has the property
$(\ast)$. Then, roughly speaking, we prove that the triple
$(A,B,W)$ is determined by $B$.
\end{rem}

\begin{thm}\label{thm:HisrepresentableinGRV}
The functor $\underline\H^\infty(e_1,\dots, e_r)$ is representable
by a closed subscheme $\H^\infty(e_1,\dots, e_r)$ of $\grv$.
\end{thm}

\begin{proof}
Recall that Theorem~\ref{thm:injectivityKricheverH} states that
$\underline\H^\infty(e_1,\dots, e_r)$ is a subfunctor of
$\underline\M^\infty(r)$. Our task consists of proving that it is
a closed subfunctor, since $\M^\infty(r)$ is closed in $\grv$
(\cite{MP}).

Theorem~\ref{thm:traceBinB} reduces the proof to checking that for
${\mathcal B}\in\M^\infty(r)(S)$ the condition $\tr({\mathcal
B})\subseteq {\mathcal B}$ is fulfilled on a closed subscheme of
$S$. Recall that such a condition is closed because ${\mathcal
B}\in\grv(S)$ and $\tr({\mathcal B})$ is quasi-coherent.
\end{proof}

\begin{rem}
Note that the points of $\H^\infty(e_1,\dots, e_r)$ corresponding
to coverings where the curve upstairs is not connected are
precisely the points of the intersection $\H^\infty(e_1,\dots,
e_r)\cap \gr^{\text{dec}}(V)$ (see
Definition~\ref{defn:decomposablegrv}).
\end{rem}

\begin{thm}\label{thm:hurwitzbilinear}
Let $B\in\M^{\infty}(r)\subset \gr^m(V)$ ($m\neq \frac12(r-n)$) be
a closed point. Let $u_1,\dots, u_r$ be integer numbers defined by
$v_m=z_1^{u_1}\dots z_r^{u_r}$.

Then, $B\in \H^{\infty}(e_1,\dots,e_r)$ if and only if the
following ``bilinear identities'' are satisfied:
    {\small
    $$\res_{z=0}\left(\Big(
    \sum_{j=1}^r\sum_{i=1}^{e_j}\frac{\psi_{u,B}^{(j)}(\xi_j^i
    z^{1/e_j},t)}{(\xi_j^i z^{1/e_j})^{\delta_{uj}-u_j}} \Big)\cdot\Big(
    \sum_{j=1}^r\sum_{i=1}^{e_j}\frac{\psi_{v,B}^{*(j)}(\xi_j^i
    z^{1/e_j},s)}{(\xi_j^i z^{1/e_j})^{u_j-1+\delta_{vj}}}
    \Big)\right)\frac{dz}z\,=\,0$$}
for all $1\leq u,v\leq r$.
\end{thm}

\begin{proof}
Let us observe that from the bilinear identity given in
Theorem~\ref{thm:bilinearidentity}, the condition $\tr(B)\subseteq
B$ is equivalent to the conditions:
    {\small $$
    T_2\Big(\tr\big(\frac{v_m
    \psi_{u,B}(z_\centerdot,t)}{(1,\dots,z_u,\dots,1)}\big),
    \frac{v_{r-n-m} \psi_{v,B}^*(z_\centerdot,s)}{(1,\dots,z_v,\dots,1)}\Big)\,=\,0
    $$}
for all $1\leq u,v\leq r$. Recalling the definition of $T_2$ and
the explicit expression of the trace of $V$ as a
$\C((z))$-algebra, one concludes.
\end{proof}

\begin{rem}
Let us observe that the bilinear identities does not characterize
the points of $\H^{\infty}(e_1,\dots,e_r)$ in $\gr^m(V)$; in fact,
a point of $\H^{\infty}(e_1,\dots,e_r)$ is characterized by these
bilinear equations and the equations characterizing
$\M^\infty(r)\subset \gr^m(V)$ (see~\cite{MP}) which are not a
hierarchy of soliton equations. This is clarified in
Theorem~\ref{thm:diffequationsofHinGRV}.
\end{rem}

\begin{thm}\label{thm:hurwitzequations}
A closed point $B\in\M^\infty(r)$ ($B\notin
\gr^{\frac12(r-n)}(V)$)  is a point of
$\H^{\infty}(e_1,\dots,e_r)$ if and only if its $\tau$-function
fulfills the following set of equations:
    {\small
    $$\begin{aligned} &
    \sum_{\substack{1\leq j\leq r \\1\leq k\leq r}}
    \Big(\sum_{\substack{1\leq i\leq e_j \\1\leq l\leq e_k}} \sum\;
    \xi_j^{i(u_j-\delta_{uj}-\alpha_1+\beta_1)} D_{\lambda_j,\alpha_1}(-\wtilde
    \partial_{t^{(j)}}) p_{\beta_1}(\wtilde\partial_{t^{(j)}})
    \check{D}^j_{\lambda_\centerdot}(-\wtilde\partial_t) \; \cdot
    \\
    & \qquad\cdot \xi_k^{l(1-u_k-\delta_{kv}-\alpha_2+\beta_2)}
    D_{\mu_k,\alpha_2}(\wtilde
    \partial_{s^{(k)}}) p_{\beta_2}(-\wtilde\partial_{s^{(k)}})
    \check{D}^k_{\mu_\centerdot}(\wtilde\partial_s) \Big)
    \tau_{B_{uj}}(t)\cdot\tau_{B_{kv}}(s) \,=\,0
    \end{aligned}$$}
for all Young diagrams $\lambda_1,\mu_1,\dots,\lambda_r,\mu_r$ and
$1\leq u,v\leq r$. Here $D_{\lambda,\alpha}$ and
$\check{D}^j_{\lambda_\centerdot}$ are differential operators
whose explicit expression will be given in the proof. The third
sum runs over the set of $4$-tuples
$\{\alpha_1,\beta_1,\alpha_2,\beta_2\}$ of non-negative integers
such that
$\frac{u_j-\delta_{uj}-\alpha_1+\beta_1}{e_j}+\frac{1-u_k-\delta_{kv}-\alpha_2+\beta_2}{e_k}=0$,
and $u_i$ are given by $v_m=z_1^{u_1}\dots z_r^{u_r}$.
\end{thm}

\begin{proof}
The proof is similar to that of Theorem~5.4 of~\cite{MP} and is
based on that given by Fay in~\cite{Fay}. To begin, let us state
some results to be used.

First, let $\chi_\lambda$ be the Schur polynomial corresponding to
$\lambda$, $t$ (resp. $s$) be the set of variables
$(t_1,t_2,\dots)$ (resp. $(s_1,s_2,\dots)$), and
$\wtilde\partial_{t}$ denotes $(\frac{\partial}{\partial
t_1},\frac{\partial}{2\partial t_2},\dots)$. From the fact
$\chi_\lambda(\wtilde\partial_t)\chi_\mu(t)\vert_{t=0}=\delta_{\lambda,\mu}$
we have that a function $f(t)\in \C\{\{t\}\}$ admits an expansion
of the type:
    {\small
    $$
    f(t)\,=\, \big(\sum_{\lambda}
    \chi_\lambda(t)\chi_\lambda(\wtilde\partial_{s})
    f(s)\big)\vert_{s=0}
    $$}
where $\lambda$ runs over the set of Young diagrams.

Second, recall Pieri's formula (\cite{Mc}, formula~I.5.16):
{\small
$$
\chi_\lambda(t) p_m(t)\,=\, \sum_{\mu-\lambda=(m)} \chi_\mu(t)$$}
where the condition $\mu-\lambda=(m)$ means that the skew diagram
$\mu-\lambda$ is a horizontal $m$-strip and $p_m(t)$ is defined by
the identity $p_m(t)=\chi_{(m)}(t)$  or, equivalently, by: {\small
\beq
 \label{eq:exponentialexpassion}
 \sum_{k\geq0} p_k(t)z^k\,=\,
\exp(\sum_{k\geq 1}t_k z^k) \end{equation}  }

 Third, the following computation is useful: {\small
    $$
    \begin{aligned}
    \chi_\lambda(\wtilde\partial_t)\big(p_m(t)f(t)\big)\vert_{t=0}
    \,&=\, \chi_\lambda(\wtilde\partial_t)\Big(\sum_\mu
    p_m(t)\chi_\mu(t)\chi_\mu(\wtilde\partial_{s})f(s)\vert_{s=0}\Big)\vert_{t=0}
    \,=\\
    &=\,
    \chi_\lambda(\wtilde\partial_t)\Big(\sum_\mu\sum_{\gamma-\mu=(m)}
    \chi_\gamma(t)\chi_\mu(\wtilde\partial_{s})f(s)\vert_{s=0}\Big)\vert_{t=0}
    \,=\\
    &=\, \sum_{\lambda-\mu=(m)}\chi_\mu
    (\wtilde\partial_{s})f(s)\vert_{s=0}
    \end{aligned}
    $$}
and let us define the operator
$D_{\lambda,m}(\wtilde\partial_{t})$ by the following identity:
    $$
    D_{\lambda,m}(\wtilde\partial_{t})f(t)\,:=\,
    \sum_{\lambda-\mu=(m)}\chi_\mu
    (\wtilde\partial_{t})f(t)\vert_{t=0}
    $$

Finally, observe that~\ref{defn:expressionBAj} allows us to write
down the following explicit expression for
$\psi_{u,B}^{(j)}(z_j,t)$:
    {\small \beq \label{eq:BAjexpanssion}
    \begin{aligned}
    \psi_{u,B}^{(j)}(z_j,t)\,& =\, \exp\big(-\sum_{i\geq
    1}\frac{t_i^{(j)}}{z_j^i}\big) \frac{\tau_{B_{uj}}(t+[z_j])}{\tau_B(t)} \,=\\
    & =\, \Big(\sum_{i\geq 0}\frac{p_i(-t^{(j)})}{z_j^i}\Big)
    \frac{\big(\sum_{i\geq 0}
    p_i(\wtilde\partial_{t^{(j)}})z_j^i\big)\tau_{B_{uj}}(t)} {\tau_B(t)}
    \end{aligned}
    \end{equation}}
since $\exp\big(\sum_{i\geq 1} z_j^i
\wtilde\partial_{t_i^{(j)}}\big)\tau_{B_{uj}}(t) =
\tau_{B_{uj}}(t+[z_j])$. Similarly, one has that:
    {\small
    $$
    \psi_{v,B}^{*(k)}(z_k,t)\,=\, \Big(\sum_{i\geq
    0}\frac{p_i(t^{(k)})}{z_k^i}\Big) \frac{\big(\sum_{i\geq 0}
    p_i(-\wtilde\partial_{t^{(k)}})z_k^i\big)\tau_{B_{kv}}(t)} {\tau_B(t)}
    $$}

We are now ready to prove the statement. The bilinear identity of
Theorem~\ref{thm:hurwitzbilinear} says that the coefficient of
$z^{-1}$ of a certain function vanishes. Note that the coefficient
of $z_j^{m}$ in $\psi_{u,B}^{(j)}(\xi_j^i z^{1/e_j},t)$ is equal
to $\xi_j^{im}$ times the coefficient of $z_j^{m}$ in
$\psi_{u,B}^{(j)}(z_j,t)$. From
formulae~\ref{eq:exponentialexpassion} and~\ref{eq:BAjexpanssion},
the residue condition reads:
    {\small $$
    \sum_{\substack{ 1\leq j\leq r \\1\leq k\leq r}}\;
    \sum_{\substack{ 1\leq i\leq e_j \\1\leq l\leq e_k}}
    \sum\Big(
    \frac{p_{\alpha_1}(-t^{(j)}) p_{\beta_1}(\wtilde\partial_{t^{(j)}})}
    {\xi_j^{i(\delta_{uj}-u_j+\alpha_1-\beta_1)}}
    \tau_{B_{uj}}(t) \cdot
    \frac{p_{\alpha_2}(s^{(k)}) p_{\beta_2}(-\wtilde\partial_{s^{(k)}})}
    {\xi_k^{l(u_k-1+\delta_{kv}+\alpha_2-\beta_2)}}
    \tau_{B_{kv}}(s)\Big) \,=\,0
    $$}
where the third sum runs over the set of $4$-tuples
$\{\alpha_1,\beta_1,\alpha_2,\beta_2\}$ of non-negative integers
such that
$\frac{u_j-\delta_{uj}-\alpha_1+\beta_1}{e_j}+\frac{1-u_k-\delta_{kv}-\alpha_2+\beta_2}{e_k}=0$.

This is a function, $F$, on $2r$ sets of variables; namely,
$t=(t^{(1)},\dots , t^{(r)})$ and $s=(s^{(1)},\dots,s^{(r)})$. Its
vanishing is equivalent to the vanishing of:
    {\small
    $$\Big(
    \prod_{\substack{ 1\leq a \leq r \\ 1 \leq b \leq r}}
    \chi_{\lambda_a}(-\wtilde
    \partial_{t^{(a)}})\chi_{\mu_b}(\wtilde
    \partial_{s^{(b)}})\Big) F\vert_{ t=s=0}\,=\,0
    $$}
for all Young diagrams $\lambda_1,\mu_1,\dots,\lambda_r,\mu_r$.

Using the facts discussed at the beginning of the proof, we arrive
at the following identity:
    {\small
    $$\begin{aligned} &
    \sum_{\substack{ 1\leq j\leq r \\1\leq k\leq r}}\;
    \sum_{\substack{ 1\leq i\leq e_j \\1\leq l\leq e_k}} \sum\Big(
    \xi_j^{i(u_j-\delta_{uj}-\alpha_1+\beta_1)} D_{\lambda_j,\alpha_1}(-\wtilde
    \partial_{t^{(j)}}) p_{\beta_1}(\wtilde\partial_{t^{(j)}})
    \check{D}^j_{\lambda_\centerdot}(-\wtilde\partial_t) \tau_{B_{uj}}(t)\;
    \cdot
    \\
    & \qquad\qquad\cdot \xi_k^{l(1-u_k-\delta_{kv}-\alpha_2+\beta_2)}
    D_{\mu_k,\alpha_2}(\wtilde
    \partial_{s^{(k)}}) p_{\beta_2}(-\wtilde\partial_{s^{(k)}})
    \check{D}^k_{\mu_\centerdot}(\wtilde\partial_s) \tau_{B_{kv}}(s) \Big)
    \,=\,0
    \end{aligned}$$}
where $\check{D}^j_{\lambda_\centerdot}(\wtilde\partial_t)$ is the
operator $\prod_{a\neq j} \chi_{\lambda_a}(\wtilde
\partial_{t^{(a)}})\vert_{ t^{(a)}=0}$.
\end{proof}

\begin{rem}
The technique used in the above proof allows one to translate some
of our previous results, such as
Theorem~\ref{thm:characterizationgrred}, into a set a differential
equations.
\end{rem}

\begin{thm}\label{thm:diffequationsofHinGRV}
Let $U$ be a closed point of $ \gr^m(V)$ ($m\neq \frac12(r-n)$)
and let $\{u_1,\dots, u_r\}$ be integer numbers defined by
$v_m=z_1^{u_1}\cdot\dots\cdot z_r^{u_r}$.

Then, $U$ is a point of $\H^\infty(e_1,\dots,e_r)$ if and only if
its Baker-Akhiezer function satisfies the following differential
equations:
\begin{itemize}
\item
the equations of Theorem~\ref{thm:hurwitzequations};
\item
the equations:
    {\small
    $$ \begin{aligned} &
    \sum_{j=1}^r   \sum\Big( D_{\lambda_j,\alpha_1}(-\wtilde
    \partial_{t^{(j)}}) p_{\beta_1}(\wtilde\partial_{t^{(j)}})
    \check{D}^j_{\lambda_\centerdot}(-\wtilde\partial_t)
    \cdot
    D_{\mu_j,\alpha_2}(-\wtilde
    \partial_{{t'}^{(j)}}) p_{\beta_2}(\wtilde\partial_{{t'}^{(j)}})
    \check{D}^j_{\mu_\centerdot}(-\wtilde\partial_{t'})
    \cdot
    \\
    & \qquad\qquad  \cdot D_{\nu_j,\alpha_3}(\wtilde
    \partial_{{t''}^{(j)}}) p_{\beta_3}(-\wtilde\partial_{{t''}^{(j)}})
    \check{D}^j_{\nu_\centerdot}(\wtilde\partial_{t''}) \Big)
    \tau_{U_{uj}}(t) \tau_{U_{vj}}(t') \tau_{U_{jw}}(t'') \,=\,0
    \end{aligned}
    $$}
for all $1\leq u,v,w\leq r$, all Young diagrams $\lambda,\mu,\nu$
and all $t,t',t''$ (the inner sum runs over the 6-tuples
$\{\alpha_1,\beta_1,\alpha_2,\beta_2,\alpha_3,\beta_3\}$ such that
$-\alpha_1+\beta_1-\alpha_2+\beta_2-\alpha_3+\beta_3=\delta_{uj}+\delta_{vj}+\delta_{wj}-u_j$)
.
\item
the equations:
    {\small
    $$
    \sum_{j=1}^r \; \sum\Big(
    D_{\lambda_j,\alpha}(\wtilde
    \partial_{t^{(j)}}) p_{\beta}(-\wtilde\partial_{t^{(j)}})
    \check{D}^j_{\lambda_\centerdot}(\wtilde\partial_t) \Big)
    \tau_{U_{ju}}(t) \,=\,0
    $$}
for all $1\leq u\leq r$, all Young diagrams $\lambda$ and all $t$
(the inner sum runs over the pairs $(\alpha,\beta)$ such that
$-\alpha+\beta=u_j+\delta_{uj}-1$).
\end{itemize}
\end{thm}

\begin{proof}
It will suffice to check that the second and third sets of
differential equations are the equations characterizing
$\M^\infty(r)$ in $\grv$. From  Theorem~\ref{thm:minftyr}, we know
that $\M^\infty(r)$ consists of those $U\in\grv$ such that $U\cdot
U\subseteq U$ and $\C\subset U$.

Theorem~\ref{thm:BAgeneragrmV} implies that these two conditions
are equivalent to:
    {\small $$
    \res_{z=0} \tr\big(\frac{v_m \psi_{u,U}(z_\centerdot,t)\psi_{v,U}(z_\centerdot,t')\psi_{w,U}^*(z_\centerdot,t'')}
    {(1,\dots, z_u,\dots, z_v,\dots,z_w,\dots,1)}
    \big)\frac{dz}{z}
    \,=\,0
    $$}
and
    {\small $$
    \res_{z=0}
    \tr\big(\frac{v_{r-n-m} \psi_{u,U}^*(z_\centerdot,t)}{(1,\dots, z_u,\dots,1)}
    \big) dz\,=\,0
    $$}
respectively. Proceeding
as in the previous proof, one concludes.
\end{proof}


\section{Curves with prescribed involutive series}
\label{sec:curveswithseries}

Theorems~\ref{thm:traceBinB}, \ref{thm:hurwitzbilinear},
\ref{thm:hurwitzequations} completely characterize  those
algebraic cur\-ves that are a finite covering of another curve
with invariants $(n; e_1,\dots, e_r)$. We shall now  apply these
results to characterize the existence of algebraic $1$-dimensional
series on a curve with prescribed numerical invariants. However,
an explicit resolution of this problem would require us to compute
the Baker-Akhiezer function of $\tr(B)$ as a point of
$\gr(\C((z)))$ for $B\in\H^{\infty}(e_1,\dots,e_r)$.

Let $Y$ be a smooth algebraic curve of genus $g$ over $\C$. An
involutive algebraic series of genus $g_0$ and degree $n$ over $Y$
is the algebraic series $\gamma_n^1$ defined by a finite morphism:
$$\pi\colon Y\,\longrightarrow\, X$$
where $X$ is a smooth algebraic curve of genus $g_0$:
$$\gamma_n^1\,=\, \{\pi^{-1}(x)\,\vert\, x\in X\}\,\subset \,
S^nY$$ or, equivalently:
$$\gamma_n^1\,\equiv \, \Gamma_\pi \,\hookrightarrow \, X\times
Y$$ ($\Gamma_\pi$ being the graph of the morphism $\pi$).

If $X$ is the projective line ($g_0=0$), then the algebraic series
defined by $\pi$ is a linear series, $g^1_n$, of degree $n$.

For instance, a curve with a linear series $g_2^1$ is a
hyperelliptic curve; a curve of genus $g>3$ with a linear series
$g_3^1$ is a trigonal curve; a curve with an algebraic series
$\gamma_2^1$ of genus $g_0>0$ is called a $g_0$-hyperelliptic
curve (\cite{Accola}).

The simplest case is the moduli space of curves of genus $g$ with
a linear series $g_n^1$ (this problem is trivial for big enough
$n$). Let us denote by $\underline\H^\infty(g,0;e_1,\dots,e_r)$
the subfunctor of $\underline\H^\infty(e_1,\dots, e_r)$ consisting
of coverings of the type:
$$\pi\colon Y \,\longrightarrow\, {\mathbb P}_1$$
where  $Y$ is of arithmetic genus $g$, $x\in {\mathbb P}_1$ and
$\pi^{-1}(x)=e_1y_1+\dots+ e_r y_r$ (with $e_1+\dots +e_r=n$).

In other words, the set
$\underline\H^\infty(g,0;e_1,\dots,e_r)(\C)$ is the set of
 curves of genus $g$ with a linear series $g_n^1$ and a
divisor $D\in g_n^1$ of the type $D=e_1y_1+\dots+ e_r y_r$.

\begin{thm}\label{thm:hg0representable}
The functor $\underline\H^\infty(g,0;e_1,\dots,e_r)$ is
representable by a closed subscheme,
$\H^\infty(g,0;e_1,\dots,e_r)$, of $\H^\infty(e_1,\dots,e_r)$.
\end{thm}

\begin{proof}
The condition that the fibres of the family of curves $Y\to S$
have arithmetic genus $g$ means that
$\underline\H^\infty(g,0;e_1,\dots,e_r)$ lies inside the connected
component $\gr^{1-g}(V)$, which is a closed subscheme of $\grv$.

The second condition, namely that $X_s={\mathbb P}_1$ for all
$s\in S$, is equivalent to saying that $\tr(K(Y))$ lies in the
connected component of index $1$, $\gr^1(\C((z)))$. This is also a
closed condition.
\end{proof}

In particular, if we set $e_1=\dots=e_r=1$ and $r=n$, then the
moduli space $\H^\infty(g,0;1,\dots,1)$ represents all curves of
genus $g$ with a linear series $g_n^1$ and parameters at the
marked points.

Let $\underline\H(g,0;1,\dots,1)$ be the Hurwitz functor
classifying the set of data $(Y,y_1,\dots,y_n)$ of coverings $Y\to
{\mathbb P}_1$ with a distinguished fibre of pairwise different
points $\{y_1,\dots,y_n\}$.

Note that there is a canonical morphism: {\small
\beq\label{eq:G+principalbundle} \Phi\colon
\underline\H^\infty(g,0;1,\dots,1)\,\longrightarrow\,
\underline\H(g,0;1,\dots,1)
\end{equation}}
that forgets the formal parameters.

\begin{thm}
The functor $\underline\H(g,0;1,\dots,1)$ is representable by a
closed subscheme, $\H(g,0;1,\dots,1)$, of
$\H^\infty(g,0;1,\dots,1)$.
\end{thm}

\begin{proof}
Let us define a morphism:
$$ \sigma\colon  \underline\H(g,0;1,\dots,1)\,\longrightarrow \,
\underline\H^\infty(g,0;1,\dots,1)$$ as follows:
$\sigma(Y,y_1,\dots,y_n)$ is the unique ${\mathcal
Y}\in\underline\H^\infty(g,0;1,\dots,1)$ on the fibre
$\Phi^{-1}(Y,y_1,\dots,y_n)$ such that $K({\mathcal Y})\cap
\C((z))=\C[z^{-1}]$.

Geometrically, this construction corresponds to choosing the set
of data $(Y,{\mathbb P}_1,\pi,x,\bar y,t_x,t_{\bar y})$ (with
$\bar y=\{y_1\,\dots,y_r\}$) such that $({\mathbb P}_1,x,t_x)$
satisfies: {\small $$ t_x\big(H^0({\mathbb P}_1-x,\o_{{\mathbb
P}_1})\big)= \C[z^{-1}]\in\gr(\C((z)))$$}

Since $\sigma$ is injective, it is enough to show that
$\underline\H(g,0;1,\dots,1)$ is a closed subfunctor of
$\underline\H^\infty(g,0;1,\dots,1)$ (via the morphism $\sigma$).
Since an $S$-valued point, $\mathcal Y$, of
$\underline\H^\infty(g,0;1,\dots,1)$, belongs to
$\underline\H(g,0;1,\dots,1)$ precisely when $\tr(K({\mathcal
Y}))\subseteq \o_S[z^{-1}]$ and since this is a closed condition
(because $\o_S[z^{-1}]$ is a point of $\gr(\C((z)))$), the theorem
is proved.
\end{proof}

\begin{rem}
Let us note that the morphisms $\Phi$ and $\sigma$ define
morphisms between the schemes $\H(g,0;1,\dots,1)$ and $
\H^\infty(g,0;1,\dots,1)$ such that $\sigma$ is a canonical
section of $\Phi$.

Furthermore,  the group $G:=\underline\aut(\C((z)))$ (see
\cite{MP2} for its definition and properties) acts on
$\gr(\C((z)))$ and on $\grv$ because, in our case, the
$\C((z))$-algebra $V$ is $\C((z))\times\dots\times \C((z))$. If
$G_+$ is the subgroup of $G$ representing ``coordinate changes''
(\cite{MP2}), then the morphism~\ref{eq:G+principalbundle} is a
$G_+$-principal bundle, that is: {\small \beq
\H(g,0;1,\dots,1)\,=\, \H^\infty(g,0;1,\dots,1)/G_+
\end{equation}}
\end{rem}

\begin{thm}
The equations defining the subscheme $\H(g,0;1,\dots,1)$ as a
subscheme of $\H^\infty(g,0;1,\dots,1)$ (via the section $\sigma$)
are as follows:
    $$\res_{z=0}\Big(z^{-i}\cdot
    \tr\big(\frac{v_{1-g} \psi_{u,B}({z_\centerdot},t)}
    {(1,\dots,z_u,\dots,1)}\big)\Big)dz
    \,=\,0
    \qquad i\geq 2\; , \;  1\leq u\leq r
    $$
\end{thm}

\begin{proof}
Note that a point $B\in \H^\infty(g,0;1,\dots,1)$ lies in
$\H(g,0;1,\dots,1)$ if and only if $\tr B= \C[z^{-1}]$. We can
write this condition in terms of Baker-Akhiezer functions.

The equality is equivalent to saying that $z^{-i}\in (\tr
B)^\perp$ for all $i\geq 2$. And the claim follows.
\end{proof}

In order to study the moduli space of curves of genus $g$ with an
involutive algebraic series $\gamma_n^1$ of genus $g_0$ and degree
$n$, one has to consider the subfunctor
$\underline\H^\infty(g,g_0;1,\dots, 1)$ of
$\underline\H^\infty(1,\dots,1)$ consisting of those coverings
$Y\to X$ where $Y$ has arithmetic genus $g$ and $X$ has arithmetic
genus $g_0$. Analogously to Theorem~\ref{thm:hg0representable},
one proves that this subfunctor is representable by a closed
subscheme of $\H^\infty(1,\dots,1)$, which will be denoted by
$\H^\infty(g,g_0;1,\dots, 1)$, since we know that the following
two conditions:
\begin{itemize}
\item $B\in \gr^{1-g}(V)$,
\item $\tr(B)\in\gr^{1-g_0}(\C((z)))$.
\end{itemize}
are closed ($B$ belongs to $\H^\infty(1,\dots,1)$).

If we assume that $B\in\gr^{1-g}(V)$, the second condition can be
translated, in some particular cases, into a set of differential
equations. We shall study this problem elsewhere and shall obtain,
for those cases, explicit characterizations (for instance,
characterizations of $g_0$-hyper\-elliptic curves).

\begin{rem}
One can consider the Hurwitz functor parametrizing families of
$n$-sheeted coverings $Y\to X$ of smooth connected curves with $Y$
of genus $g$ and $X$ of genus $g_0$. Let us denote it by
$\underline\H_n(g,g_0)$. Let ${\mathcal W}$ denote the open
subscheme of $\H^\infty(g,g_0;1,\dots, 1)$ whose points
corresponds to geometric data with both curves  smooth and
connected. We therefore have a canonical forgetful morphism:
$$\Psi\colon {\mathcal W}\,\longrightarrow\, \underline\H_n(g,g_0)$$
Note that the fibre of $\Psi$ over a covering, $Y\to X$, is
isomorphic to the complementary in $X$ of the set of ramification
points.
%
\end{rem}

\section{Final Remarks}

It was pointed out by Li and Mulase~\cite{LiMulaseHitchin} that a
Zariski open subset of the moduli space of Higgs pairs over a
curve can be embedded into a quotient Grassmannian and that the
restriction of the $n$-component KP-flow is precisely the
Hamiltonian flow of the Hitchin system.

In our setting, a local analogue for the Hitchin system appears
naturally when we try to solve the following question: how do the
algebro-geometric objects introduced here ($\tau$-functions, Baker
functions, etc.) depend on the structure of $\C((z))$-algebra of
$V$?. Let us discuss  this question briefly.

For a monic polynomial of degree $n$, $P(T)\in \C[[z]][T]$, with
coefficients in $\C[[z]]$, let us define the finite
$\C((z))$-algebra: $$V_P\,:= \, \C((z))/(P(T))$$
 and a rank $n$ free
$\C[[z]]$-module $(V_P)_+:= \C[[z]]/(P(T))$. Then, $\w C_P:=\sf
(V_P)_+$ is called the formal spectral cover of polynomial $P$.

Let us denote by $\A_{\infty}^n$ the infinite dimensional affine
space representing $\C[[z]]\times\overset{n}\dots\times \C[[z]]$
(see~\cite{MP},~\S3.A). For each point $s=(s_1,\dots,
s_n)\in\A_{\infty}^n$, let us denote by $P_s$ the polynomial
$T^n-s_1 T^{n-1}+\dots+(-1)^n s_n$ and, for the sake of
simplicity, let us write $V_s=V_{P_s}$, $V_{s+}=(V_{P_s})_+$ and
$\w C_s=\w C_{P_s}$.

One can define a family of infinite Grassmannians parametrized by
the space $\A_{\infty}^n$: \beq \label{eq:grfibration} {\mathcal
G}r \,\overset{\pi}\longrightarrow\, \A_{\infty}^n \end{equation}
 such
that the fibre of $s$, $\pi^{-1}(s)$, is the infinite Grassmannian
of the couple $(V_s,V_{s+})$.

There is an open dense subscheme $U\subset \A_{\infty}^n$ such
that for each $s\in U$ the formal spectral cover $\w C_s$ is
smooth. The fibres of $\pi$ over $U$ correspond to  the
Grassmannians of $(V_s,V_{s+})$ where $V_s$ is a separable
$\C((z))$-algebra, which have been studied in~\S\S1-4 of this
paper.

Note that there is a natural representation of $\End_{\C[[z]]}V_+$
in the Lie algebra of vector fields over $\grv$:
$$ \End_{\C[[z]]}V_+ \,\hookrightarrow\, \End_{\C((z))}V \,\overset{\Psi}\longrightarrow\, T \grv$$
where  the fibre of $\Psi$ at the point $U$ is given by: {\small
$$
\begin{aligned}
\Psi_U \colon \End_{\C((z))}V & \,\longrightarrow\, T_U
\grv\,=\,\hom(U,V/U) \\
\varphi & \,\longmapsto \, \Psi_U(\varphi):= \big(
U\hookrightarrow V\overset{\varphi}\to V\to V/U\big)
\end{aligned}$$}

Using this representation, a notion of a local Higgs pair can be
introduced and the corresponding moduli space can be studied. The
local analogue of the Hitchin fibration is then related to the
fibration~(\ref{eq:grfibration}). Further, the different
$\underline n$-KP hierarchies can be interpreted as the flows
defining the fibres of the local Hitchin map. We hope to study all
these aspects elsewhere.

\begin{rem}
Let us denote by $V_{s_1}$ and $V_{s_2}$ two different
$\C((z))$-algebra structures on $V$ and let us assume that these
structures are determined by two partitions of $n$, $\underline
n_1=(e_1^1,\dots, e_{r_1}^1)$ and $\underline n_2=(e_1^2,\dots,
e_{r_2}^2)$. In~\S9 of~\cite{AdamsBergvelt}, the following
question is stated: given a point $U$ coming from two geometric
data: {\small
$$(Y^i,X^i,\pi^i, x^i,\bar y^i, t_{x^i}, t_{\bar
y^i})\,\in\, \H^\infty(e_1^1,\dots, e_{r_i}^1)\,\subset
\,\grv\qquad i=1,2
$$} is there any relation between the curves $Y^1$ and $Y^2$?. The
answer to this question is that in general $Y^1$ and $Y^2$ do not
have to be isomorphic (even if both curves are irreducible).

The following example, related to the constructions
of~\S\ref{sec:curveswithseries}, is instructive for understanding
the above question. Let $\pi^i\colon Y^i\to {\mathbb P}^1$
($i=1,2$) be two different finite coverings of degree $n$ such
that:
$$\pi^1_*\o_{Y^1}\,\simeq\, \pi^2_*\o_{Y^2}$$
Obviously, we can find non-isomorphic smooth curves $Y^1$ and
$Y^2$ fulfilling this condition, since the set of rank $n$ locally
free sheaves on ${\mathbb P}^1$ is discrete (they are of the form
$\o_{{\mathbb P}^1}(a_1)\oplus\dots\oplus \o_{{\mathbb
P}^1}(a_n)$) while the set of degree $n$ coverings is not.
Choosing a point $x\in{\mathbb P}^1$ and a formal trivialization
at that point, $t_x$, one observes that $(Y^1,{\mathbb
P}^1,x,t_x)$ and $(Y^2,{\mathbb P}^1,x,t_x)$ define the same point
of $\grv$.
\end{rem}


\end{document}